\documentclass[10pt]{article}
\usepackage{amsmath,amssymb}

\newcommand{\bra}{\langle}
\newcommand{\ket}{\rangle}
\newcommand{\ep}{\hfill {$\Box$}}

\newtheorem{thm}{Theorem}[section]
\newtheorem{cor}{Corollary}[section]
\newtheorem{defin}{Definition}[section]
\newtheorem{lem}{Lemma}[section]
\newtheorem{prop}{Proposition}[section]

\title{Commutation Relations for Unitary Operators}
\author{\small{M.A. Astaburuaga, O. Bourget, V.H. Cort\'es \footnote{Supported by the Grants Fondecyt 1120786, ECOS-Conicyt C10E10}}\\
\small{Facultad de Matem\'aticas, Pontificia Universidad Cat\'olica de Chile,}\\
\small{Av. Vicu\~na Mackenna 4860, Macul, Santiago, Chile}\\
\small{E-mail: bourget@mat.puc.cl}\\
\small{phone: (56 2) 2354 4509}\\
}
\date{}

\begin{document}
\maketitle

\begin{abstract}
Let $U$ be a unitary operator defined on some infinite-dimensional complex Hilbert space ${\cal H}$. Under some suitable regularity assumptions, it is known that a local positive commutation relation between $U$ and an auxiliary self-adjoint operator $A$ defined on ${\cal H}$ allows to prove that the spectrum of $U$ has no singular continuous spectrum and a finite point spectrum, at least locally. We show that these conclusions still hold under weak regularity hypotheses and without any gap condition. As an application, we study the spectral properties of the Floquet operator associated to some perturbations of the quantum harmonic oscillator under resonant AC-Stark potential.
\end{abstract}


\section{Introduction}

The spectral analysis of unitary operators defined on Hilbert spaces is a natural tool in the study of the long-time behavior of periodic time-dependent quantum systems \cite{ev}. It also appears in the theory of orthogonal polynomials \cite{S1}, \cite{S2} and the study of classical dynamical systems e.g \cite{N1}, \cite{N2}.

The commutation relations satisfied by an operator may be relevant to determine its spectral properties. This approach has been developped to a large extent for self-adjoint operators to analyze either its discrete spectrum \cite{hs2}, \cite{lp}, \cite{ah}, \cite{stubbe}, \cite{hs1}, or its essential component by means of some positive commutator methods \cite{p}, \cite{mo}, \cite{jmp}, \cite{abmg}, \cite{bm}, \cite{skib}, \cite{gg}, \cite{ggm}, \cite{gj}. The development of these methods within the spectral theory of unitary operators has been historically delayed, although this gap is now partly filled regarding the development of the positive commutator theory \cite{p}, \cite{abcf}, \cite{frt}.

In this manuscript, we show that the traditional conclusions of the Mourre theory for unitary operators still hold under weak regularity conditions and without any extra gap condition. This is synthesized by Theorem \ref{thm0}. The strategy followed for the proof is a transposition of the exposition given in \cite{abmg} to our unitary setting. The proof is intrinsically based on the unitary functional calculus. Then, this abstract result is applied to  derive new results concerning local perturbations of the Floquet operator associated to a quantum harmonic oscillator under a resonant AC-Stark potential.

The manuscript is structured as follows. The abstract result is presented in Section 2, the example treated in Section 3. The proof of Theorem \ref{thm0} is postponed to Section 4. It is almost entirely based on the unitary functional calculus. However, we have tried to emphasize the analogies between the self-adjoint and unitary cases by using the notations of \cite{abmg}. Some auxiliary results and technicalities have been gathered in Sections 5, 6, 7.\\

\noindent {\bf Notations:} Let us fix some notations adopted throughout this paper. Our unitary operator is defined on some fixed infinite-dimensional Hilbert space ${\cal H}$ on ${\mathbb C}$. The resolvent set of a closed operator $B$ on ${\cal H}$ is denoted by $\rho(B)$ and its spectrum by: $\sigma(B) \equiv {\mathbb C}\setminus \rho(B)$. The open unit disk and the unit circle are denoted by ${\mathbb D}$ and $\partial {\mathbb D}={\mathbb S}$ respectively. The one-dimensional torus is denoted by ${\mathbb T}$. The positive constants independent of the relevant parameters of the problem are generically denoted by $c$ or $C$. If $A$ is a self-adjoint operator defined on ${\cal H}$ with domain ${\cal D}(A)$, we use the japanese bracket notation: $\bra A\ket = \sqrt{(A^2+1)}$. Lastly, for any function $\Phi$ on ${\mathbb S}$ is associated in a unique manner to the function $\phi$ defined on ${\mathbb T}$ by:  $\phi(\theta)=\Phi(e^{i\theta})$, for all $\theta \in {\mathbb T}$. If $U$ is a unitary operator defined on ${\cal H}$ and if its spectral family is denoted by $(E_{\Delta})_{\Delta \in {\cal B}({\mathbb T})}$, where ${\cal B}({\mathbb T})$ stands for the family of Borel sets of ${\mathbb T}$, we will have that:
\begin{equation*}
\Phi(U)=\int_{\mathbb T}\phi(\theta)dE(\theta)=\int_{\mathbb T}\Phi(e^{i\theta})dE(\theta)\enspace .
\end{equation*}
We will identify frequently the spectrum of $U$ and its component (which are subsets of ${\mathbb S}$) with the corresponding support of the spectral measure, which lies in ${\mathbb T}$.

\section{An Abstract Result}

In this section, we introduce the main abstract result of this manuscript i.e. Theorem \ref{thm0}. The core of its development relies on the existence of a self-adjoint operator $A$, densely defined on ${\cal H}$ (the conjugate operator), which respect to which our unitary operator $U$ satisfies some suitable regularity conditions. We start by describing them.

\begin{defin}\label{c1} Let $B\in {\cal B}({\cal H})$ and $A$ a self-adjoint operator defined on ${\cal H}$ with domain ${\cal D}(A)$. The operator $B$ is of class $C^1$ with respect to $A$ (or shortly $B\in C^1(A)$), if there exists a dense linear subspace ${\cal S}$ of ${\cal H}$, ${\cal S}\subset {\cal D}(A)$, such that the sesquilinear form $F$, defined
\[F(\varphi,\phi):=\langle A\varphi,B\phi\rangle - \langle \varphi,BA\phi\rangle\]
for any $(\varphi,\phi)\in {\cal S}\times {\cal S}$, extends continuously to a bounded form on ${\cal H}\times {\cal H}$. The bounded linear operator associated to the extension of $F$ is denoted by $\mathrm{ad}_A (B) =[A,B]$.
\end{defin}

\begin{defin}\label{ck} Let $k\in {\mathbb N}$, $B\in {\cal B}({\cal H})$ and $A$ a self-adjoint operator defined on ${\cal H}$ with domain ${\cal D}(A)$. The operator $B$ is of class $C^k$ with respect to $A$ (or shortly $B\in C^k(A)$), if there exists a dense linear subspace ${\cal S}$ of ${\cal H}$, ${\cal S}\subset {\cal D}(A)$, such that:
\begin{itemize}
\item $B\in C^{k-1}(A)$
\item the sesquilinear form $F$, defined by: $F(\varphi,\phi):=\langle A\varphi, \mathrm{ad}_A^{k-1}(B) \phi\rangle - \langle \varphi, \mathrm{ad}_A^{k-1}(B) A\phi\rangle $, for any $(\varphi,\phi)\in {\cal S}\times {\cal S}$, extends continuously to a bounded form on ${\cal H}\times {\cal H}$.
\end{itemize}
The bounded linear operator associated to the extension of $F$ is denoted by $\mathrm{ad}_A (\mathrm{ad}_A^{k-1}(B)) =\mathrm{ad}_A^{k}(B)$. If $B$ belongs to $C^k(A)$ for any $k\in {\mathbb N}$, we say that $B\in C^{\infty}(A)$. 
\end{defin}

Actually, the notation takes its origin in the fact that a bounded linear operator $B$ belongs to $C^k(A)$ if and only if the strongly continuous application $t \mapsto e^{itA}B e^{-itA}$ with values in ${\cal B}({\cal H})$ is strongly $C^k$ on ${\mathbb R}$. We refer to Section 5 or \cite{abmg} for more details. We shall write naturally: $C^0(A)={\cal B}({\cal H})$ and $\mathrm{ad}_A^0 B=B$. This alternative point of view allows to consider intemediate scales of regularity leading to the definition of the classes ${\cal C}^{s,p}(A)$ (see Definition \ref{csp}, \cite{abmg} Chapter 5 or \cite{sah}). Without going into the details for the moment, we say that $B\in {\cal C}^{1,1}(A)$ if:
\begin{equation*}
\int_0^1 \|e^{iA\tau}B e^{-iA\tau}+e^{-iA\tau}B e^{iA\tau} -2B\|\,\frac{d\tau}{|\tau|^2} < \infty \enspace .
\end{equation*}
In particular, $C^2(A)\subset {\cal C}^{1,1}(A)\subset C^1(A)$.\\

\noindent{\bf Remark:} ${\cal S}$ can be equivalently chosen as ${\cal D}(A)$ in Definitions \ref{c1} and \ref{ck}. Indeed, assume that two densely defined sesquilinear forms coincide on some dense subdomain ${\cal S}\times {\cal S} \subset {\cal H}\times {\cal H}$. If one of them extends continuously to a bounded form on ${\cal H}\times {\cal H}$, the other also extends continuously to a bounded form on ${\cal H}\times {\cal H}$, and both extensions coincide. \\

If $U$ is a unitary operator defined on some Hilbert space ${\cal H}$, $U\in C^1(A)$ if and only if $U^*\in C^1(A)$. In particular, $U({\cal D}(A))$ and $U^*({\cal D}(A))$ are subsets of ${\cal D}(A)$, which implies that: $U({\cal D}(A))=U^*({\cal D}(A))={\cal D}(A)$. These considerations motivates the following equivalence, proved in Section 5:
\begin{lem}\label{equiv1} Let $U$ be a unitary operator defined on ${\cal H}$. Then, the following assertions are equivalent:
\begin{itemize}
\item[(a)] $U\in C^1(A)$.
\item[(b)] $U^*\in C^1(A)$.
\item[(c)] There exists a dense linear subspace ${\cal S}_1$ of ${\cal H}$ such that $U {\cal S}_1= {\cal S}_1$, ${\cal S}_1 \subset {\cal D}(A)$ and the sesquilinear form $F_1:{\cal S}_1\times {\cal S}_1\rightarrow {\mathbb C}$: $F_1(\varphi,\phi):=\langle U\varphi,AU\phi\rangle\ - \langle\varphi,A\phi\rangle$ extends continuously to a bounded form on ${\cal H}\times {\cal H}$. This extension is associated to a bounded operator denoted by $U^{*}AU-A$.
\item[(d)] There exists a dense linear subspace ${\cal S}_2$ of ${\cal H}$ such that $U {\cal S}_2= {\cal S}_2$, ${\cal S}_2 \subset {\cal D}(A)$ and the sesquilinear form $F_2:{\cal S}_2\times {\cal S}_2\rightarrow {\mathbb C}$: $F_2(\varphi,\phi):=\langle \varphi,A\phi\rangle - \langle U^*\varphi,AU^*\phi\rangle$ extends continuously to a bounded form on ${\cal H}\times {\cal H}$. This extension is associated to a bounded operator denoted by $A-UAU^*$.
\end{itemize}
Moreover, $U^{*}AU-A = U^* (\mathrm{ad}_A U)$, $(\mathrm{ad}_A U) U^*=A-UAU^*$.
\end{lem}
For an alternative characterization of the classes $C^k(A)$ in the unitary context, see also Lemma \ref{equiv2}. \\

When speaking about the positivity conditions we are about to introduce, we will write indifferently $U^{*}AU-A$ for $U^* (\mathrm{ad}_A U)$ and $(\mathrm{ad}_A U) U^*$ for $A-UAU^*$ (in the sense of Lemma \ref{equiv1}):
\begin{defin}\label{propagating} Let $A$ be a self-adjoint operator with domain ${\cal D}(A) \subset {\cal H}$ and $U$ a unitary operator which belongs to $C^1(A)$. Then, we say that given $\Theta \in {\cal B}({\mathbb T}),$
\begin{itemize}
\item {\bf $P_w({\mathbb T})$}: $U$ is weakly propagating with respect to $A$ if $U^{*}AU-A > 0$ (i.e non-negative and injective) 
\item {\bf $P({\Theta})$}: $U$ is propagating with respect to the observable $A$ on ${\Theta}$ or on the arc $e^{i\Theta}$ if there exist $c >0$ and a compact operator $K$ such that: $E_{\Theta} (U^{*}AU-A) E_{\Theta} \geq c E_{\Theta} + K$
\item {\bf $P_s(\Theta)$}: $U$ is strictly propagating with respect to the observable $A$ on ${\Theta}$ or on the arc $e^{i\Theta}$ if there exist $c >0$ such that: $E_{\Theta} (U^{*}AU-A) E_{\Theta} \geq c E_{\Theta}$.
\end{itemize}
\end{defin}
We have clearly that: ${\bf P_s({\mathbb T})}\Rightarrow {\bf P_w({\mathbb T})}$. Sometimes, we write that the operator $U$ is (strictly) propagating for $A$ at a point $\theta$ of the torus ${\mathbb T}$, when there exists an open neighbourhood $\Theta_{\theta}$ of $\theta$ such that $U$ is (strictly) propagating for $A$ on ${\Theta}_{\theta}$. Following \cite{abcf}, this is equivalent to claim that there exist a smoothed characteristic function $\phi$ supported in $\Theta_{\theta}$, which takes value 1 on a neighbourhood of $\theta$ and a positive constant $c$ such that: 
\begin{equation*}
\Phi(U)\left(U^*AU-A\right)\Phi(U) \geq c \Phi(U)^2\enspace.
\end{equation*}

\noindent{\bf Remark:} Since the spectral projectors associated to $U$ commute with $U$ and $U^*$, the positivity conditions presented in Definition \ref{propagating} can be equivalently described writing $A-UAU^*$ in place of $U^*AU-A$. This remark will be used without any further comment.

We can formulate now a first spectral result:
\begin{thm}\label{thm2} Assume that the unitary operator $U$ is weakly propagating with respect to the self-adjoint operator $A$. Then, $\sigma_{pp}(U)=\emptyset$.
\end{thm}
This result was proven under a somewhat different form in \cite{abcf2}. We give a straightforward proof based on the Virial Theorem in Paragraph 4.1.

However, strenghtening the regularity hypotheses, we can derive more precise informations on the spectral properties of $U$. The following result was exposed initially in \cite{abcf} and then proven under weaker hypotheses in \cite{frt}:
\begin{thm}\label{thm} Let $\Theta$ be an open subinterval of ${\mathbb T}$. Assume that $U$ is propagating with respect to $A$ on $\Theta$ and that $U^{*}AU-A \in C^1(A)$. Then,
\begin{itemize}
\item $U$ has a finite number of eigenvalues in $e^{i\Theta}$. Each of these eigenvalues has a finite multiplicity.
\item For any $\theta \in \Theta\setminus \sigma_{pp}(U)$, there exists $\delta>0$ such that:
\begin{equation*}
\sup_{|z|\neq 1, \arg z \in (\theta-\delta, \theta+ \delta)} \|\bra A \ket^{-1}(1-zU^*)^{-1}\bra A \ket^{-1} \| < \infty \enspace .
\end{equation*}
\item The spectrum of $U$ has no singular continuous component in $e^{i\Theta}$.
\end{itemize}
\end{thm}
The control of the point spectrum was again obtained by a suitable version of the Virial Theorem and only uses the fact that $U$ is propagating with respect to $A$ on $\Theta$ (see Paragraph 4.1). \\

The extension of Theorem \ref{thm} in \cite{frt} was obtained by associating the existing theory for self-adjoint operator with the Cayley transform. However, this approach is still limited by an extra spectral gap condition, which may be a handicap for the applications (see Section 3). In this manuscript, we show that this limitation can be overcome with Theorem \ref{thm0}. 

Following the notations of \cite{trieb}, \cite{abmg} Chapter 2 and denoting by ${\cal K}$ the interpolation space $({\cal D}(A),{\cal H})_{1/2,1}$ (continuously and densely embedded in ${\cal H}$), we have that:
\begin{thm}\label{thm0} Let $\Theta$ be an open subset of ${\mathbb T}$. Assume $U$ is propagating with respect to $A$ on $e^{i\Theta}$ and that $U\in {\cal C}^{1,1}(A)$. Then, 
\begin{itemize}
\item[(i)] $U$ has a finite number of eigenvalues in $e^{i\Theta}$. Each of these eigenvalues has a finite multiplicity. Denoting by $e^{i\Theta'}=e^{i\Theta}\setminus \sigma_{\text{pp}}(U)$, we have that:
\item[(ii)] Given any $(\psi,\varphi)\in {\cal K}^2$, the holomorphic functions defined on ${\mathbb D}$ by $z\mapsto \bra \varphi, (1-zU^*)^{-1}\psi\ket$ and $z\mapsto \bra \varphi, (1-\bar{z}^{-1}U^*)^{-1}\psi\ket$ extend continuously to ${\mathbb D}\cup e^{i\Theta'}$.
\item[(iii)] The spectrum of $U$ has no singular continuous component in $e^{i\Theta}$.
\end{itemize}
If $U$ is strictly propagating with respect to $A$ on $e^{i\Theta}$, then Statements (i) and (iii) can be replaced by: $U$ is purely absolutely continuous on $e^{i\Theta}$.
\end{thm}

\noindent {\bf Remark:} Statement (ii) can be reformulated equivalently as follows: the limits
$\lim_{r\rightarrow 1^-}\bra \varphi, (1-re^{i\theta}U^*)^{-1}\psi\ket$ and $\lim_{r\rightarrow 1^-} \bra \varphi, (1-r^{-1}e^{i\theta}U^*)^{-1}\psi\ket$ exist uniformly when $\theta$ belongs to any compact subset $K \subset \Theta'$. Denoting these limits respectively by $F^{\pm}(\theta)$ we deduce that the maps $\theta \mapsto F^{\pm}(\theta)$ are continuous on $\Theta'$ and the corresponding weak$^*$ limit operators, $w^*$-$\lim_{r\rightarrow 1^-}(1-re^{i\theta}U^*)^{-1}$ and $w^*$-$\lim_{r\rightarrow 1^-} (1-r^{-1}e^{i\theta}U^*)^{-1}$, are well-defined bounded operators in ${\cal B}({\cal K},{\cal K}^*)$. \\

\noindent {\bf Remark:} Since $C^2(A)\subset C^{1,1}(A)$, Theorem \ref{thm0} includes Theorem \ref{thm}. \\

The next section is dedicated to the example. The proof of Theorem \ref{thm0} is postponed to Section 4.

\section{On Resonant AC-Stark perturbations for the Harmonic Oscillator}

Time-dependent perturbations of the quantum harmonic oscillator have been a regular subject of interest \cite{hls}, \cite{ev}, \cite{comb}, \cite{gy}, \cite{wmw}. When it is submitted to an AC-Stark potential, the Floquet operator of the system, which is explicit, undergoes a spectral transition between the resonant and non-resonant regimes \cite{ev}. The stability of these spectral properties under perturbations has been studied in the non-resonant regime \cite{wmw} and partly in the resonant regime \cite{gy}. Both approaches are based on the Floquet Hamiltonian formalism. Our aim is to investigate some perturbations in the resonant regime, which are not covered by \cite{gy}.

\subsection{Preliminaries and Main Results}

We recall the main features of the model briefly. For more details, we refer the reader to \cite{ev}. The Hamiltonian of the Harmonic oscillator (with unit mass) is defined on $L^2({\mathbb R})$ by:
$$
H_{\omega} = \frac{p^2}{2} + \frac{1}{2}\, \omega^2 x^2
$$
where $p = -i \partial_x$. We also write $\omega_0 T = 2\pi$. Let $E$ be a real-valued continuous periodic function with period $T$, $T>0$ and $(U_0 (t,s))_{(s,t)\in {\mathbb R}^2}$ be the unitary propagator associated to the AC-Stark Hamiltonian $H_0 (t) = H_{\omega} + E(t)\, x$. The propagator is explicit. In particular, for all $t\in {\mathbb R},$
\begin{equation}
U_0(t,0) = e^{-i \varphi_1(t) x}e^{i \varphi_2(t) p/\omega}e^{-iH_{\omega}t -i\psi(t)} \label{free}
\end{equation}
where $\varphi_1(t) = \int_0^t E(\tau) \cos(\omega (\tau - t)) \, d\tau$, $\varphi_2(t) =  -\int_0^t E(\tau) \sin(\omega (\tau - t ))\, d\tau$ and $\psi(t) =-\frac{1}{2}\int_0^t (\varphi_1(\tau)^2-\varphi_2(\tau)^2)\, d\tau$. For simplicity, we denote $U_0(t,0)$ by $U_0(t)$ in the following. The evolution of the observables $x$ and $p$ under this propagator are also explicit (understood on a suitable domain like the space of the Schwartz functions ${\cal S}({\mathbb R})$):
\begin{eqnarray}\label{man2}
U_0(t)^* p \, U_0(t) & = &  -x \omega \sin(\omega t) + p \cos(\omega t) + \varphi_1(t) \label{freeEnss}\\
U_0(t)^* x \, U_0(t) & = &  x \cos(\omega t) + \frac{p}{\omega} \sin(\omega t) - \frac{1}{\omega}  \varphi_2(t)\nonumber
\end{eqnarray}
These identities allows us to deduce the spectral properties of the Floquet operator $U_0(T)$. Indeed, paraphrasing \cite{ev} Proposition 4.1, if $\omega_0 \neq \omega$, the functions $\varphi_1$ and $\varphi_2$ are bounded on ${\mathbb R}$, meaning that the orbits $\{U_0(t)\psi; t\in {\mathbb R}\}$ are precompact for all vector $\psi \in {\cal H}$. As a consequence,
\begin{prop} If $\omega_0 \neq \omega$, the Floquet operator $U_0(T)$ is pure point.
\end{prop}

If $\omega_0 =\omega$, situation coined as resonant, the spectral properties of $U_0(T)$ can be drastically modified. From now, we write $A_1=\varphi_1 (T)^{-1}p$ if $\varphi_1 (T)\neq 0$ and $A_2=-\omega \varphi_2 (T)^{-1}x$ if $\varphi_2 (T)\neq 0$.
\begin{prop}\label{freeho} Let $\omega_0 =\omega$. Then, the following hold:
\begin{itemize}
\item[(a)] If $\varphi_1 (T)=\varphi_2 (T)=0$, then $U_0(T)$ is pure point.
\item[(b)] If either $\varphi_1 (T)\neq 0$ or $\varphi_2 (T)\neq 0$, then $U_0(T)$ has purely absolutely continuous spectrum and $\sigma(U_0(T))={\mathbb S}$. Specifically, noting that $U_0(T){\cal S}({\mathbb R})={\cal S}({\mathbb R}),$
\begin{itemize}
\item If $\varphi_1 (T)\neq 0$, then the operator $U_0(T)^* A_1 U_0(T) - A_1$ defined via its sesquilinear form on ${\cal S}({\mathbb R})\times {\cal S}({\mathbb R})$ can be extended uniquely as a bounded operator on $L^2({\mathbb R})$ and
\begin{equation}\label{freeEnss2}
U_0(T)^* A_1 U_0(T) - A_1 = I \enspace ,
\end{equation}
\item If $\varphi_2 (T)\neq 0$, then the operator $U_0(T)^* A_2 U_0(T) - A_2$ defined via its sesquilinear form on ${\cal S}({\mathbb R})\times {\cal S}({\mathbb R})$ can be extended uniquely as a bounded operator on $L^2({\mathbb R})$ and
\begin{equation*}
U_0(T)^* A_2 U_0(T) - A_2 = I \enspace .
\end{equation*}
\end{itemize}
\item[(c)] $U_0(T)$ belongs to $C^{\infty}(A_1)\cap C^{\infty}(A_2)$ or equivalently to $C^{\infty}(p)\cap C^{\infty}(x)$.
\end{itemize}
\end{prop}
\noindent{\bf Proof:} Statement (a) follows from the explicit form of the propagator (\ref{free}). Statements (b) and (c) follow from relations (\ref{freeEnss}), Theorem \ref{thm} and the considerations on the Heisenberg couples developped in Section 6. The fact that $U_0(T)$ belongs to $C^{\infty}(A_1)\cap C^{\infty}(A_2)$ or equivalently to $C^{\infty}(p)\cap C^{\infty}(x)$ can also be derived directly. Actually, for all $k\in {\mathbb N},$
\begin{eqnarray}
\mathrm{ad}_p^k U_0(T) &=& \varphi_1 (T)^k U_0(T)\label{freeEnss3}\\
\mathrm{ad}_x^k U_0(T) &=& (-\frac{\varphi_2 (T)}{\omega})^k U_0(T) \nonumber \enspace.
\end{eqnarray}
\ep

Let us mention that the proof of the absolute continuity is obtained in \cite{ev} and \cite{gy} for a specific example ($E(t)=\sin \omega t$) although the procedure can be easily generalized in the latter case. Our approach using positive commutators allows us not only to extend this result to a wider class of electric field $E$ but also to consider some perturbations of the model as shown in this section.

As mentioned at the beginning, we show that in the resonant regime some of the spectral property remains if the Hamiltonian $H_0(\cdot)$ is perturbed. Let $V$ denote the multiplication operator by the real-valued function $V(\cdot)$ on $L^2({\mathbb R})$, $V(\cdot) \in L^{\infty}({\mathbb R})$ and define the perturbed time-dependent Hamiltonian, $H(\cdot)$ by: $H(t)= H_0(t) + V$. If the propagator $(U(t,s))$ associated to $H(\cdot)$ exists, then it satisfies for all $(s,t)\in {\mathbb R}^2$, $U(t,s) = U_0(t,s) \Omega(t,s)$ where $\Omega (t,s)$ is defined in the strong sense by:
\begin{equation}\label{dyson1}
\Omega(t,s)-I = -i \int_s^t U_0(\tau,s) V U_0^*(\tau,s) \Omega(\tau,s)\, d\tau \enspace,
\end{equation}
or equivalently by its Dyson expansion (see e.g. \cite{rs} Theorem X.69):
\begin{equation}\label{dyson2}
\Omega(t,s) = I + \sum_{j=1}^{\infty} (-i)^j \int_s^t \ldots \int_s^{\tau_{j-1}} U_0(\tau_1,s) V U_0^*(\tau_1,s) \dots U_0(\tau_j,s) V U_0^*(\tau_j,s) \, d\tau_j \ldots d\tau_1 \enspace.
\end{equation}
As before, we denote $\Omega(t):=\Omega(t,0)$ and $U(t):=U(t,0)$.

The spectral properties of $U(T)$ are described by Theorem \ref{a}:
\begin{thm}\label{a} Let $\omega_0 =\omega$ and assume that $\varphi_j (T)\neq 0$ for some $j\in \{1,2\}$. Given a real-valued function $V(\cdot)$ such that:
\begin{equation*}
\int_0^1 \left( \int_{\mathbb R}|V(x-t) +V(x+t)-2V(x) |\, dx \right) \,\frac{dt}{t^2} < \infty \enspace,
\end{equation*}
consider the Floquet operator $U(T)$ defined by (\ref{dyson1}). Then,
\begin{itemize}
\item[(a)] If $\partial_x V(x)$ vanishes when $|x|$ tends to infinity, then there is no singular continuous component in the spectrum of $U(T)$. Moreover, its point subspace has finite dimension.
\item[(b)] if $T \|\partial_x V(\cdot)\|_{\infty} < |\varphi_1 (T)|$ (resp. $2\pi \|\partial_x V(\cdot)\|_{\infty} < |\varphi_2 (T)|$), then the spectrum of $U(T)$ is purely absolutely continuous.
\end{itemize}
Moreover, given any $(\psi,\varphi)\in {\cal K}^2$, the holomorphic functions defined on ${\mathbb D}$ by $z\mapsto \bra \varphi, (1-zU^*)^{-1}\psi\ket$ and $z\mapsto \bra \varphi, (1-\bar{z}^{-1}U^*)^{-1}\psi\ket$ extend continuously to $\overline{{\mathbb D}}\setminus \sigma_{pp}(U)$.
\end{thm}

Statement (b) was proven in \cite{gy} for smooth and mildly unbounded potentials $V$. To our knowledge, statement (a) is new.

The proofs of Theorem \ref{a} is an application of Theorem \ref{thm0}. The purpose of the next paragraphs is to relate the hypotheses of Theorem \ref{a} with the framework of Theorem \ref{thm0}. For simplicity, we assume in the following that $\varphi_1(T)\neq 0$. The other case can be treated similarly.

\subsection{Regularity issues}

In this paragraph, we relate the regularity properties of $V$ to the regularity properties of $U(T)$ with respect to $A_1$ (or equivalently to $p$):

\begin{lem}\label{lemEnss0} Let $n \in {\mathbb N}$ and $V(\cdot)$ be a real-valued function in $C^n({\mathbb R})$ such that for all $k\in \{0,\ldots,n\}$, $\partial_x^k V(\cdot) \in L^\infty (\mathbb{R})$. Then, for all $t\in {\mathbb R}$, $U_0^*(t) V U_0(t)$ belongs to $C^n(p)$. Precisely, for all $t\in {\mathbb R}$ and all $k\in \{0,\ldots,n\},$
\begin{eqnarray*}
\mathrm{ad}^k_p(U_0^*(t)V U_0(t)) &=& (-i)^k\, cos^k(\omega t) U^*_0(t) \partial_x^k V U_0(t) \\
\| \mathrm{ad}_p^{k}(U_0^*(t) V U_0(t)) \| &\leq & \|\partial_x^k V\|_\infty \enspace.
\end{eqnarray*}
\end{lem}
\noindent{\bf Proof:} Fix $t\in {\mathbb R}$. Using sesquilinear forms on ${\cal S}({\mathbb R})\times {\cal S}({\mathbb R})$, we have that for any function $W(\cdot)$ in $C^1({\mathbb R})$, $W\in C^1(U_0(t) p U_0^*(t))$ and
\begin{equation*}
\mathrm{ad}_p (U_0^*(t)WU_0(t)) = U_0^*(t) \left( \mathrm{ad}_{U_0(t) p U_0^*(t)} W \right) U_0(t)= (-i)\, cos(\omega t) U_0^*(t) \partial_x W U_0(t) \enspace,
\end{equation*}
due to identity \ref{man2}. The proof follows by induction on $k$, $k\in \{0,\ldots\, n\}$. For $k=0$, the conclusion is straightforward since we are dealing with bounded operators. Once the conclusion established for some $k\in \{0,\ldots\, n-1\}$, we observe that $\partial_x^k V(\cdot)\in C^1({\mathbb R})$. It follows from the previous considerations that: $\partial_x^k V\in C^1(U_0(t) p U_0^*(t))$. Using the induction hypothesis, we have that:
\begin{equation*}
\mathrm{ad}_p^{k+1}(U_0^*(t)VU_0(t))=\mathrm{ad}_p \left( \mathrm{ad}_p^{k}(U_0^*(t) V U_0(t))\right) = (-i)^k\, cos^k(\omega t) \mathrm{ad}_p U^*_0(t) \partial_x^k V U_0(t) \enspace,
\end{equation*}
which leads to the conclusion. \ep

\begin{lem}\label{lemEnss0-1} Let $n \in {\mathbb N}$ and $V(\cdot)$ be a real-valued function in $C^n({\mathbb R})$ such that for all $k\in \{0,\ldots,n\}$, $\partial_x^k V(\cdot) \in L^\infty (\mathbb{R})$. Then, for all $t\in {\mathbb R}$, $\Omega(t)$ belongs to $C^n(p)$. In particular, for all $t\in {\mathbb R}$, $\Omega^*(t) p \Omega(t) -p$ is bounded and:
\begin{eqnarray}
\Omega^*(t) p \Omega(t) -p &=& - \int_0^t \cos (\omega \tau ) U^*(\tau ) (\partial_x V) U (\tau )\, d\tau \label{eq3}\\
\|\Omega^*(t) p \Omega(t) -p\| &\leq & |t| \| \partial_x V(\cdot)\|_\infty \label{bound1} \enspace.
\end{eqnarray}
\end{lem}
\noindent{\bf Proof:} The first part of this lemma results from Lemma \ref{lemEnss0} and Corollary \ref{dyson}. Since $\Omega(t)$ belongs to $C^1(p)$, this means that $\Omega(t) {\cal D}(p) \subset {\cal D}(p)$ in virtue of Proposition \ref{2}. Therefore, using quadratic forms on ${\cal D}(p)\times {\cal D}(p)$, we have that:
\begin{eqnarray*}
\Omega^*(t) p \Omega(t) -p & = & \int_0^t \frac{d}{d\tau} (\Omega^*(\tau) p \Omega (\tau)) \, d\tau \\
&=& - i \int_0^t  \Omega^*(\tau) [p, U_0^*(\tau) V U_0 (\tau)] \Omega(\tau) \, d\tau \\
&=& - \int_0^t \cos (\omega \tau ) \Omega^*(\tau) U_0^* (\tau) (\partial_x V) U_0 (\tau)  \Omega (\tau)\, d\tau \enspace,
\end{eqnarray*}
making use again of Lemma \ref{lemEnss0}. Since both sides define bounded operators on $L^2({\mathbb R})$, the previous identity can be extended uniquely on $L^2({\mathbb R})\times L^2({\mathbb R})$ and the conclusion is straightforward. \ep

\begin{lem}\label{lemEnss1} Let $n\in {\mathbb N}$ and $V(\cdot)$ be a real-valued function in $C^n({\mathbb R})$ such that for all $k\in \{0,\ldots,n\}$, $\partial_x^k V(\cdot) \in L^\infty (\mathbb{R})$. Then, $U(T)$ belongs to $C^n(p)$. In particular, $U^*(T) p U(T) - p$ is bounded and
\begin{equation}\label{B(T)}
U^*(T) p U(T) - p = \varphi_1(T) - \int_0^T \cos (\omega \tau ) U^* (\tau) (\partial_x V) U (\tau)\, d\tau \enspace .
\end{equation}
\end{lem}
\noindent{\bf Proof:} The first part follows from Lemma \ref{lemEnss0-1}, Proposition \ref{4} and the fact that $U_0(T)$ belongs to $C^{\infty}(p)$ (see Proposition \ref{freeho}). Using identity (\ref{freeEnss2}) we have that:
\begin{equation*}
U^*(T)pU(T)-p= \varphi_1(T)I + \Omega^*(T)p \Omega(T) -p \enspace,
\end{equation*}
which implies the result by Lemma \ref{lemEnss0-1}. \ep

As noted in Corollary \ref{equiv2}, $U(T)$ belongs to $C^n(p)$ if and only if the operator $B(T):=U^*(T)pU(T)-p$ defined by (\ref{B(T)}) belongs to $C^{n-1}(p)$.

The proof of the next result is more involved and has been postponed to Paragraph 5.5.
\begin{lem}\label{lemEnss3} If $V \in {\cal C}^{1,1}(p)$ i.e.
\begin{equation*}
\int_0^1 \left( \int_{\mathbb R}|V(x-t) +V(x+t)-2V(x) |\, dx \right) \,\frac{dt}{t^2} < \infty \enspace,
\end{equation*}
then, $U(T) \in {\cal C}^{1,1}(p)$.
\end{lem}

\subsection{Compactness issues}

In this paragraph, we relate the decay properties of $V$ or $\partial_x V$ to some compactness properties.

\begin{lem}\label{lemEnss2-0} Let $T>0$ and $W$ be a continuous function defined on ${\mathbb R}\times [0,T]$ such that:
\begin{equation*}
\lim_{|x| \to \infty} \sup_{t\in [0,T]} | W(x,t) | = 0 \enspace.
\end{equation*}
If for each $t\in [0,T]$, $W_t$ denotes the multiplication operator by the function $W(\cdot,t)$, then for all $(s,t)\in [0,T]^2$ the strongly convergent integral
\begin{equation*}
\int_s^t U_0^*(\sigma) W_{\sigma} U_0(\sigma) \, d\sigma
\end{equation*}
defines a compact operator on $L^2({\mathbb R})$.
\end{lem}
\noindent{\bf Proof:} Fix $(s,t)\in {\mathbb R}^2$ and denote
\begin{equation*}
X(t,s) = \int_s^t U_0^*(\sigma) W U_0(\sigma) \, d\sigma \enspace.
\end{equation*}
Using \cite{ev} Lemma 5.4, the assertion is true for $W$ in $C_0^{\infty}({\mathbb R}\times [0,T])$.  Now, let us take
a sequence $W_n \in C_0^{\infty}({\mathbb R}\times [0,T])$ such that $\sup_{(x,t) \in {\mathbb R}\times [0,T]}|W_n(x,t)-W(x,t)| \to 0$, as $n$ tends to infinity. We denote $X_n(t,s) = \int_s^t U_0^*(\sigma) (W_n)_{\sigma} U_0(\sigma) \, d\sigma$. So, $X_n(t,s)$ is a compact operator which converges in norm to $X(t,s)$. Indeed,
\begin{eqnarray*}
X_n(t,s) -  X(t,s) &=& \int_s^t U_0^*(\sigma) ((W_n)_{\sigma}-W_{\sigma}) U_0(\sigma) \, d\sigma \\
\|X_n(t,s) - X(t,s) \| &\leq & |t-s| \, \|W_n - W\|_{\infty} \enspace,
\end{eqnarray*}
which proves that $X(t,s)$ is compact. \ep

It follows from Lemma \ref{lemEnss2-0} and \cite{ev} Theorem 5.2 that:
\begin{lem} Assume that $V(\cdot) \in C^0({\mathbb R})$ and that:
\begin{equation*}
\lim_{|x| \to \infty} V(x) = 0 \enspace.
\end{equation*}
Let $U$ and $U_0$ be the unitary propagators related by (\ref{dyson1}). Then, $U(t) - U_0(t)$ is compact for every $t\in {\mathbb R}$.
\end{lem}

However, we shall not use this result in our present analysis. The next result is borrowed from \cite{abm}:
\begin{lem}\label{lemEnss2} Let $T>0$ and $t\mapsto W_t$ a strongly continuous family of bounded linear operator such that for all $(s,t)\in [0,T]^2$, the operator
\begin{equation*}
\int_s^t U_0^*(\tau) W_{\tau} U_0(\tau)\, d\tau
\end{equation*}
is compact. Then, the strongly convergent integral
\begin{equation*}
\int_0^T U^*(\tau) W_{\tau} U(\tau) \, d\tau
\end{equation*}
defines also a compact operator.
\end{lem}
\noindent{\bf Proof:} Due to the Banach-Steinhaus Theorem, $\sup_{t\in [0,T]}\|W_t\|< \infty$. Let $N\in {\mathbb N}$ and define $I_j = \frac{T}{N}[j,j+1)$, $j\in \{0,\ldots,N-1\}$. On one hand, we have that:
\begin{eqnarray*}
\|U^*(t)U_0(t)-U^*(s)U_0(s)\| &=& \|U(s,t)U_0(t,s)-I\| \leq \int_0^T \|U_0^*(\tau) W_{\tau} U_0(\tau)\|\, d\tau \\
&\leq & |t-s| \sup_{t\in [0,T]}\|W_t\| \enspace ,
\end{eqnarray*}
for all $(s,t)$ in $[0,T]^2$. In particular, there exists $C>0$ such that given $j\in \{0,\ldots,N\}$ and any $(s,t,\tau)\in I_j^3$, $|t-s|\leq C N^{-1}$, $\|U^*(\tau)U_0(\tau)-U^*(s)U_0(s)\|\leq C N^{-1}$ and
\begin{equation*}
\|\int_s^t U^*(\tau) W_{\tau} U(\tau) \, d\tau - U^*(s)U_0(s) \left(\int_s^t U_0^*(\tau) W_{\tau} U_0(\tau) \, d\tau\right)U_0^*(s)U(s)\|\leq C N^{-2} \enspace.
\end{equation*}
The result follows since
\begin{equation*}
\| \int_0^T U^*(\tau) W_{\tau} U(\tau)\, d\tau - \sum_{j=0}^{N-1} U^*(\frac{jT}{N})U_0(\frac{jT}{N}) \int_{I_j} U_0^*(\tau) W_{\tau} U_0(\tau)\, d\tau U_0^*(\frac{jT}{N})U(\frac{jT}{N}) \| \leq \frac{C}{N} \enspace .
\end{equation*}
\ep

As an immediate corollary of Lemmata \ref{lemEnss0-1}, \ref{lemEnss2-0} and \ref{lemEnss2}, we have that:
\begin{cor}\label{lemEnss2-2} Assume that $V(\cdot) \in C^1({\mathbb R})$ and that:
\begin{equation*}
\lim_{|x| \to \infty} \partial_x V(x) = 0 \enspace.
\end{equation*}
Then, the bounded operator $\Omega^*(t) p \Omega(t) -p$ is compact for any $t\in {\mathbb R}$.
\end{cor}

\subsection{Proof of Theorem \ref{a}}

The hypothesis
\begin{equation*}
\int_0^1 \left( \int_{\mathbb R}|V(x-t) +V(x+t)-2V(x) |\, dx \right) \, \frac{dt}{t^2} < \infty \enspace,
\end{equation*}
means that $V\in C^{1,1}(p)$ (see Definition \ref{csp}). Therefore, $V$ belongs to $C^1(A_1)$ and we deduce that:
\begin{itemize}
\item $U(T)$ is propagating with respect to $A_1$ on ${\mathbb T}$ if
\begin{equation*}
\lim_{|x| \to \infty} \partial_x V(x) = 0 \enspace .
\end{equation*}
\item $U(T)$ is strictly propagating with respect to $A_1$ on ${\mathbb T}$ if $T \|\partial_x V(\cdot)\|_{\infty} < |\varphi_1 (T)|$.
\end{itemize}
On the other hand, Lemma \ref{lemEnss3} implies that $U(T)$ belongs also to ${\cal C}^{1,1}(A_1)$. The conclusion follows naturally from Theorem \ref{thm0}.

\section{On the proof of Theorem \ref{thm0}}

The proof of Theorem \ref{thm0} follow the lines of its self-adjoint counterpart \cite{abmg} Chapter 7. Its development is articulated on two axes:
\begin{itemize}
\item The control of the (embedded) point spectrum by means of the Virial Theorem (Paragraph 4.1)
\item The study of the continuous component of the spectrum using Mourre differential inequality strategy (Paragraph 4.2)
\end{itemize}
The proof is carried out in Paragraph 4.3.

Before starting, let us remind or fix some notations. ${\mathbb D}^*$ will stand for ${\mathbb D}-\{ 0\}$. If $\Theta$ is an open interval in ${\mathbb T}$ and $r>1$, we denote by $S_{\Theta,r}^{\pm}$ and $\Omega_{\Theta,r}^{\pm}$ the sectors
\begin{eqnarray*}
S_{\Theta,r}^+ &=& \{ z \in {\mathbb C}; \, \arg(z) \in \Theta, r^{-1}<|z|<1\}\\
S_{\Theta,r}^- &=& \{ z \in {\mathbb C}; \, \arg(z) \in \Theta, 1<|z|<r\}\\
\Omega_{\Theta,r}^+ &=& \{ z \in {\mathbb C}; \, \arg(z) \in \Theta, r^{-1}<|z|\leq 1\}\\
\Omega_{\Theta,r}^- &=& \{ z \in {\mathbb C}; \, \arg(z) \in \Theta, 1\leq |z|<r\}\enspace.
\end{eqnarray*}
The spectral measure of $U$ is denoted by $(E(\Delta))_{\Delta \in {\cal B}({\mathbb T})}$.

\subsection{The Virial Theorem and its consequences}

As mentionned at the beginning of this section, the control of the point spectrum is achieved after establishing the Virial Theorem. Its proof can be found in \cite{abcf}. We scheme out an alternative version here.
\begin{thm}\label{virial} Assume that $U\in C^1(A)$. Then, for all $\theta \in {\mathbb T},$
\begin{equation*}
E_{\{\theta\}} (U^* AU -A) E_{\{\theta\}} =0 \enspace .
\end{equation*}
In particular, if $\varphi$ is an eigenvector of $U$, $\bra \varphi, (U^* AU -A)\varphi \ket =0$.
\end{thm}
\noindent{\bf Proof:} For all $(\varphi,\psi)$ in ${\cal H}\times {\cal H},$
\begin{equation*}
\bra \varphi, (U^*AU-A) \psi \ket = \frac{i}{t}\lim_{t\rightarrow 0}\bra \varphi, (U^*e^{-iA t}U - e^{-iA t}) \psi\ket\enspace.
\end{equation*}
On the other hand, for all $(\varphi,\psi)$ in ${\cal H}\times {\cal H}$ and all $t\in {\mathbb R},$
\begin{equation*}
\bra \varphi, E_{\{\theta\}}(U^* e^{-iA t}U-e^{-iA t}) E_{\{\theta\}} \psi\ket =0 \enspace,
\end{equation*}
which implies the result. \ep

As a first consequence, we can prove Theorem \ref{thm2}.\\
\noindent{\bf Proof of Theorem \ref{thm2}:} Let $\theta \in {\mathbb T}$ such that, $E(\{\theta\})\neq 0$ and $\varphi$ an associated (non trivial) eigenvector of $U$: $E(\{\theta\})\varphi=\varphi$. Since $U$ satisfies also $P_w({\mathbb T})$, we have that: $\bra \varphi,(U^*AU-A) \varphi\ket>0$, which contradicts Theorem \ref{virial}. \ep

We also deduce that:
\begin{cor}\label{virial-2} Assume $U$ is propagating with respect to $A$ on the Borel subset $\Theta\subset {\mathbb T}$. Then, $U$ has a finite number of eigenvalues in $\Theta$. Each of these eigenvalues has finite multiplicity. 
\end{cor}
We refer to \cite{abcf} Corollary 5.1 for a proof.

If $U$ is propagating with respect to $A$ on some Borel subset $\Theta\subset {\mathbb T}$, it follows that $\Theta \cap \sigma_{\text{pp}}(U)$ is finite. Therefore, for any $\theta \in \Theta\setminus \sigma_{\text{pp}}(U)$, there exist $\delta_{\theta} >0$ and $c_{\theta}>0$ such that:
\begin{equation*}
E_{(\theta-2\delta_{\theta}, \theta+2\delta_{\theta})} (U^* AU -A) E_{(\theta-2\delta_{\theta}, \theta+2\delta_{\theta})} \geq c_{\theta} E_{(\theta-2\delta_{\theta}, \theta+2\delta_{\theta})} \enspace .
\end{equation*}
In other words, $U$ is strictly propagating at $\theta$. This motivates the development of the next section.

Let us mention another consequence of the Virial Theorem. This is an straightforward adaptation of \cite{abmg} Lemma 7.2.12. We give the proof for completeness:
\begin{lem} Let $\Theta$ be an open subset of ${\mathbb T}$, $a\in {\mathbb R}$ and $K$ be a compact operator such that:
\begin{equation}\label{e1}
E(\Theta) (U^*AU-A) E(\Theta)\geq a E(\Theta)+K \enspace .
\end{equation}
Then, for all $(\theta,\eta) \in \Theta \times (0,\infty)$, there exists $\epsilon >0$ and $F$ a finite rank operator such that: $F\leq E(\{\lambda\})$ and
\begin{equation*}
(a-\eta) (E(\Theta_{\epsilon})-F)-\eta F \leq E(\Theta_{\epsilon}) (U^*AU-A) E(\Theta_{\epsilon})
\end{equation*}
with $\Theta_{\epsilon}:=(\theta-\epsilon, \theta+\epsilon)$. If $\theta$ is not an eigenvalue: $(a-\eta) E(\Theta_{\epsilon}) \leq E(\Theta_{\epsilon}) (U^*AU-A) E(\Theta_{\epsilon})$. If $\theta$ is an eigenvalue $\min(a-\eta, \eta) E(\Theta_{\epsilon}) \leq E(\Theta_{\epsilon}) (U^*AU-A) E(\Theta_{\epsilon})$.
\end{lem}
The proof apply several times the following observation: if $(T_k)$ is a decreasing sequence of orthogonal projections which converges strongly to an operator $T$, then for any compact self-adjoint operator $C$, $\lim_{k\rightarrow \infty}\|T_k CT_k -TCT\|=0$. Therefore, given $\nu>0$, there exists $n\in {\mathbb N}$ such that: $T_n CT_n -TCT\geq -\nu I$. Multiplying on the right and the left by $T_n$ and observing that $TT_n=T=T_nT$, implies that: $T_n CT_n -TCT\geq -\nu T_n$.\\

\noindent{\bf Proof:} In order to simplify the notations let us write: $B=U^*AU-A$, $E=E(\Theta)$ and $P=E(\{\theta\})$. Moreover, for any orthogonal projection $G\leq E$, we denote: $\overline{G}=E-G$. Let $(F_k)$ be an increasing sequence of finite rank orthogonal projections such that s-$\lim_{k\rightarrow \infty} F_k =P$. Applying the previous remark, it follows that there exists a finite rank orthogonal projection such that $F\leq P$ and $\overline{F}K \overline{F} \geq \overline{P}K \overline{P}-\eta/2 \overline{F}$. Multiplying inequality (\ref{e1}) on the left and the right by $\overline{F}\leq E$ implies that:  $\overline{F}B \overline{F} \geq \overline{P}K \overline{P}+(a-\eta/2) \overline{F}$. Theorem \ref{virial} implies that: $PBP=0$. Therefore, $FBF=(P-F)BF=FB(P-B)=0$. Since $P-F=\overline{F}-\overline{P}$, $\overline{F}BF=\overline{P}BF$ and $FB\overline{F}=FB\overline{P}$. Following with algebraic considerations, $E=F+\overline{F}$ so that:
\begin{equation}\label{e2}
EBE= \overline{P}BF+\overline{F}B\overline{F}+FB\overline{P} \geq (a-\eta/2) \overline{F}+K'
\end{equation}
where the compact operator $K'=\overline{P}BF+\overline{P}K\overline{P}+FB\overline{P}$. Define now the decreasing sequence of orthogonal projections $(E_k)_{k\geq k'}$ by: $E_k=E(\Theta_{k^{-1}})$ with $k'$ chosen such that $\Theta_{k'^{-1}} \subset \Theta$ ($E_k\leq E$ and $E_k \overline{F}=E_k-F$). Multiplying inequality (\ref{e2}) on the right and left by $E_k$ implies that: $E_kBE_k\geq (a-\eta/2)(E_k-F)+E_k K' E_k$. Since s-$\lim_{k\rightarrow \infty} E_k =P$ and $PK'P=0$, it follows from the preliminar remark that we can choose $k\geq k'$ such that: $E_k K' E_k \geq -\eta/2 E_k \geq -\eta/2 (E_k-F) -\eta F$, which implies the first part of the lemma with $\epsilon =k^{-1}$. The second part is straightforward since in this case $F=0$. The last part follows observing that taking $b=\min(a-\eta, -\eta)$, $(a-\eta)(E(\Theta_{\epsilon})-F)-\eta F\geq b (E(\Theta_{\epsilon})-F)+b F = b E(\Theta_{\epsilon})$. \ep

\subsection{Differential inequalities}

What follows is an adaptation of \cite{abmg} paragraph 7.3 to our unitary formalism.

From now and until the end of this paragraph, we assume that $U$ is strictly propagating with respect to some self-adjoint operator $A$ at $\theta_0 \in {\mathbb T}$:
\begin{equation*}
E_{(\theta_0 - 2\delta, \theta_0 + 2\delta)} (A-UAU^*) E_{(\theta_0 - 2\delta, \theta_0 + 2\delta)} \geq a_1 E_{(\theta_0 - 2\delta, \theta_0 + 2\delta)} \enspace ,
\end{equation*}
for some $\delta>0$, $a_1>0$ (see also Corollary \ref{equiv1}). We also assume that $(B(\varepsilon))_{\varepsilon\in (0,\varepsilon_0]}$ and $(U_{\varepsilon})_{\varepsilon \in (0,\varepsilon_0]}$ ($\varepsilon_0>0$) are families of unitary and bounded operators on ${\cal H}$ such that:
\begin{itemize}
\item there exists $C>0$, such that for all $\varepsilon \in (0,\varepsilon_0]$, $\varepsilon \, \|B(\varepsilon)\| + \|U_{\varepsilon}-U\| \leq C \varepsilon$.
\item $\lim_{\varepsilon \rightarrow 0}\|B(\varepsilon)-(A-UAU^*)\| =0$.
\end{itemize}

Denoting $B_1:=A-UAU^*$, the following properties are straightforward:
\begin{lem}\label{lem1} There exists $C>0$ such that for all $\varepsilon \in (0,\varepsilon_0],$
\begin{gather*}
\|U_{\varepsilon}^*e^{-\varepsilon B(\varepsilon)}-U^*\| \leq C\, \varepsilon \\
\|U_{\varepsilon}^*(e^{\varepsilon B(\varepsilon)})^*-U^*\| \leq C\, \varepsilon \\
\|e^{-\varepsilon B(\varepsilon)}(e^{-\varepsilon B(\varepsilon)})^* - e^{-2\varepsilon B_1}\|\leq C\, \varepsilon \\
\|(e^{\varepsilon B(\varepsilon)})^* e^{\varepsilon B(\varepsilon)} - e^{2\varepsilon B_1}\|\leq C\, \varepsilon \enspace.
\end{gather*}
\end{lem}
Note that: $(e^{\varepsilon B(\varepsilon)})^*= e^{\varepsilon B(\varepsilon)^*}$. For $\varepsilon \in (0,\varepsilon_0]$ and $z\in \overline{\mathbb D}\setminus \{0\}$, we define:
\begin{eqnarray*}
T_{\varepsilon}^{+}(z) &=& 1- z U_{\varepsilon}^* e^{-\varepsilon B(\varepsilon)}\\
T_{\varepsilon}^{-}(z) &=& 1- \bar{z}^{-1} U_{\varepsilon}^* (e^{\varepsilon B(\varepsilon)})^*
\end{eqnarray*}

We deduce the following estimates:
\begin{lem}\label{lem3} Let $0<a_0<a_1$. There exist $\varepsilon_1 \in (0,\varepsilon_0]$, $b>0$ such that for all $z \in S^+_{(\theta_0 - \delta, \theta_0 + \delta),2}$, all $\varepsilon \in (0,\varepsilon_1]$ and all $\psi\in {\cal H},$
\begin{eqnarray*}
a_0\|\psi\|^2 & \leq & \bra \psi, B_1 \psi \ket +\frac{b}{(1-|z|)^2+\frac{16}{\pi^2}|z| \delta^2} \| T_{\varepsilon}^+(z) \psi \|^2 \\
a_0\|\psi\|^2 & \leq & \bra \psi, B_1 \psi \ket +\frac{b |z|^2}{(1-|z|)^2+\frac{16}{\pi^2}|z| \delta^2} \| T_{\varepsilon}^-(z) \psi \|^2 \enspace.
\end{eqnarray*}
\end{lem}
\noindent {\bf Proof:} Without restriction, we can assume that: $0< \delta < \pi/4$. Denote $E^{\perp}=1-E$ with $E:=E_{(\theta_0 - 2\delta, \theta_0 + 2\delta)}$. Fix $z \in S_{(\theta_0-\delta,\theta_0+\delta),0}^+$. In particular, $\bar{z}^{-1} \in S_{(\theta_0-\delta,\theta_0+\delta),0}^-$. We have that:
\begin{eqnarray*}
\|(1-zU^*)^{-1}E^{\perp}\| & \leq & d(z,\sigma(E^{\perp}UE^{\perp}))^{-1}\\
\|(1-\bar{z}^{-1}U^*)^{-1}E^{\perp}\| & \leq & d(\bar{z}^{-1},\sigma(E^{\perp}UE^{\perp}))^{-1}
\end{eqnarray*}
where $d(z,E^{\perp}UE^{\perp}) = \min d(z,e^{i(\theta_0 \pm 2\delta)})$ and $d(\bar{z}^{-1},E^{\perp}UE^{\perp}) = \min d(\bar{z}^{-1},e^{i(\theta_0 \pm 2\delta)})$.
\begin{eqnarray*}
d(z,e^{i(\theta_0 \pm 2\delta)}) &=& (1-|z|)^2 +4 |z| \sin^2 (\frac{\theta-\theta_0}{2}\mp \delta) \geq  (1-|z|)^2 +\frac{16}{\pi^2} |z| \delta^2\\
d(\bar{z}^{-1},e^{i(\theta_0 \pm 2\delta)}) &=& (1-|z|^{-1})^2 +4 |z|^{-1} \sin^2 (\frac{\theta-\theta_0}{2}\mp \delta) \geq (1-|z|^{-1})^2 +\frac{16}{\pi^2} |z|^{-1} \delta^2 \enspace.
\end{eqnarray*}
Let $\psi \in {\cal H}$. We obtain for any $\varepsilon \in (0,\varepsilon_0],$
\begin{eqnarray*}
\|E^{\perp} \psi\|^2 &=& \| (1-zU^*)^{-1}E^{\perp}(T_{\varepsilon}^{+}(z)+z(U_{\varepsilon}^*e^{-\varepsilon B(\varepsilon)}-U^*)) \psi\|^2\\
&\leq & 2\| (1-zU^*)^{-1}E^{\perp}\|^2 \left( \|T_{\varepsilon}^{+}(z) \psi\|^2 + |z|^2 \|(U_{\varepsilon}^*e^{-\varepsilon B(\varepsilon)}-U^*) \psi\|^2 \right)\\
\|E^{\perp} \psi\|^2 &=& \| (1-\bar{z}^{-1}U^*)^{-1}E^{\perp}(T_{\varepsilon}^{-}(z)+\bar{z}^{-1} (U_{\varepsilon}^*(e^{\varepsilon B(\varepsilon)})^*-U^*)) \psi\|^2\\
&\leq & 2\| (1- \bar{z}^{-1}U^*)^{-1}E^{\perp}\|^2 \left( \|T_{\varepsilon}^{-}(z) \psi\|^2 + |z|^{-2} \|(U_{\varepsilon}^* (e^{\varepsilon B(\varepsilon)})^*-U^*) \psi\|^2 \right)
\end{eqnarray*}
It follows from Lemma \ref{lem1} that there exists $C>0$ such that for all $z \in S_{(\theta_0 - \delta, \theta_0 + \delta),0}^+$, all $\varepsilon \in (0,\varepsilon_0]$ and all $\psi\in {\cal H},$
\begin{eqnarray}\label{Eperp}
\|E^{\perp} \psi\|^2 &\leq & \frac{2}{(1-|z|)^2 +\frac{16}{\pi^2} |z| \delta^2} \left( \|T_{\varepsilon}^{+}(z) \psi\|^2 + C |z|^2 \varepsilon^2 \|\psi\|^2\right) \\
\|E^{\perp} \psi\|^2 &\leq & \frac{2}{(1-|z|^{-1})^2 +\frac{16}{\pi^2} |z|^{-1} \delta^2} \left( \|T_{\varepsilon}^{-}(z) \psi\|^2 + C |z|^{-2} \varepsilon^2 \|\psi\|^2\right) \enspace .
\end{eqnarray}
On the other hand, we have that: $a_1 \|\psi \|^2 \leq a_1 \|E^{\perp} \psi \|^2 +\bra E\psi, B_1 E\psi \ket$
with $\bra E\psi, B_1 E\psi \ket = \bra \psi, B_1 \psi \ket - \bra E^{\perp}\psi, B_1 E^{\perp}\psi \ket - 2\Re (\bra E\psi, B_1 E^{\perp}\psi \ket)$. 
Given any $(\nu,\varepsilon) \in (0,1]\times (0,\varepsilon_0],$
\begin{equation*}
2 \|E \psi \| \|B_1 E^{\perp}\psi\|  \leq \nu \|E \psi \|^2 + \nu^{-1} \| B_1 E^{\perp} \psi \|^2 \leq \nu \|\psi \|^2 + \nu^{-1} \| B_1\|^2 \| E^{\perp} \psi \|^2
\end{equation*}
entailing
\begin{equation*}
\bra E\psi, B_1 E\psi \ket \leq \bra \psi, B_1 \psi \ket + \| B_1\| \| E^{\perp} \psi \|^2 + \nu \|\psi \|^2 + \nu^{-1} \| B_1\|^2 \| E^{\perp} \psi \|^2
\end{equation*}
and
\begin{equation*}
a_1 \|\psi \|^2 \leq \bra \psi, B_1 \psi \ket + \nu \|\psi \|^2 + \left[a_0 + \| B_1\|+\nu^{-1} \| B_1\|^2\right] \| E^{\perp} \psi \|^2 \enspace .
\end{equation*}
Combining with inequalities (\ref{Eperp}), we deduce that for all $\varepsilon \in (0,\varepsilon_0],$
\begin{align*}
[a_1 - \nu & - \frac{2 C |z|^2 \varepsilon^2}{(1-|z|)^2 +\frac{16}{\pi^2} |z| \delta^2}\left(a_1 + \| B_1\|+\nu^{-1} \| B_1\|^2\right) ] \|\psi \|^2 \\
& \leq \bra \psi, B_1 \psi \ket + \frac{2}{(1-|z|)^2 +\frac{16}{\pi^2} |z| \delta^2}\left(a_1 + \| B_1\|+\nu^{-1} \| B_1\|^2\right) \|T_{\varepsilon}^{+}(z) \psi\|^2 \enspace ,\\
[a_1 - \nu & - \frac{2 C \varepsilon^2}{(1-|z|)^2 +\frac{16}{\pi^2} |z| \delta^2}\left(a_1 + \| B_1\|+\nu^{-1} \| B_1\|^2\right) ] \|\psi \|^2 \\
& \leq \bra \psi, B_1 \psi \ket + \frac{2}{(1-|z|^{-1})^2 +\frac{16}{\pi^2} |z|^{-1} \delta^2}\left(a_1 + \| B_1\|+\nu^{-1} \| B_1\|^2\right) \|T_{\varepsilon}^{+}(z) \psi\|^2 \enspace .
\end{align*}
Fix $\nu=(a_1-a_0)/4$. There exists $\varepsilon_1 \in (0,\varepsilon_0]$, such that for all $\varepsilon \in (0,\varepsilon_1]$ and $z \in S_{(\theta_0 - \delta, \theta_0 + \delta),2}^+,$
\begin{equation*}
\frac{2 C \varepsilon^2}{(1-|z|)^2 +\frac{16}{\pi^2} |z| \delta^2}\left(a_1 + \| B_1\|+\nu^{-1} \| B_1\|^2\right) \leq C (\frac{\pi \varepsilon}{2 \delta})^2 \left(a_1 + \| B_1\|+\nu^{-1} \| B_1\|^2\right) \leq \frac{a_1-a_0}{4}
\end{equation*}
and the former estimates rewrite
\begin{eqnarray*}
a_0 \|\psi \|^2 & \leq & \bra \psi, B_1 \psi \ket + \frac{b}{(1-|z|)^2 +\frac{16}{\pi^2} |z| \delta^2} \|T_{\varepsilon}^{+}(z) \psi\|^2\\
a_0 \|\psi \|^2 & \leq & \bra \psi, B_1 \psi \ket + \frac{b}{(1-|z|^{-1})^2 +\frac{16}{\pi^2} |z|^{-1} \delta^2} \|T_{\varepsilon}^{+}(z) \psi\|^2 \enspace ,
\end{eqnarray*}
with $b=2(a_1 + \| B_1\|+4(a_1-a_0)^{-1} \| B_1\|^2)$. \ep

Our next step deals with the invertibility of the family of bounded operators $(T_{\varepsilon}^{\pm}(z))$. Before, let us make a couple of observations.
\begin{lem}\label{lem2} Let $A$ be a bounded invertible operator on ${\cal H}$. Then, $(1-A)$ is invertible if and only if $(1-(A^{-1})^*)$ is invertible. In this case,
\begin{equation*}
\Re ((1+A)(1-A)^{-1}) = 2 \Re ((1-A)^{-1}) -1= (1-A)^{-1} - (1-(A^{-1})^*)^{-1}\enspace.
\end{equation*}
\end{lem}
\noindent{\bf Proof:} The first remark follows from the fact that a bounded operator $B$ is invertible if and only if its adjoint $B^*$ is invertible ($(B^*)^{-1}=(B^{-1})^*$) and identity: $(1-A) = - A(1-A^{-1})$. Rewriting, 
\begin{equation*}
(1+A)(1-A)^{-1} = (1-A)^{-1}+(A-1+1)(1-A)^{-1}=2 (1-A)^{-1} -1 \enspace,
\end{equation*}
gives the first identity. The second part follows, if one notes that: $B(B-1)^{-1}=1-(1-B)^{-1}$ with $B$ in place of $A$ or $(A^{-1})^*$. \ep
\\

Our second observation is that the functions $h_1$ and $h_2$ defined on ${\mathbb R}$ by: $h_1(0)=h_2(0)=2$ and
\begin{eqnarray*}
h_1(x) &=& x^{-1} (1-e^{-2x})\\
h_2(x) &=& x^{-1} (e^{2x}-1)
\end{eqnarray*}
for $x\neq 0$, are homeomorphisms from ${\mathbb R}$ onto $(0,\infty)$, respectively monotone decreasing and monotone increasing. Since for any $\varepsilon \in (0,\varepsilon_0]$, $1-e^{-2 \varepsilon B_1}=\varepsilon B_1 h_1(\varepsilon B_1)$ and $e^{2 \varepsilon B_1}-1=\varepsilon B_1 h_2(\varepsilon B_1)$, we have that: $c_1 \varepsilon B_1 \leq 1-e^{-2 \varepsilon B_1}$ and $c_2 \varepsilon B_1 \leq e^{2 \varepsilon B_1} -1$ for some positive constant $c_1$ and $c_2$.

Both remarks are used in the following lemma:
\begin{lem}\label{lem4} The linear operators $T_{\varepsilon}^{\pm}(z)$ are invertible in ${\cal B}({\cal H})$, provided $(\varepsilon,z) \in (0,\varepsilon_2]\times \Omega_{(\theta_0 - \delta, \theta_0 + \delta),2}^+$ or $(\varepsilon,z) \in [0,\varepsilon_2]\times S_{(\theta_0 - \delta, \theta_0 + \delta),2}^+$ for some $\varepsilon_2 \in (0,\varepsilon_1]$. Denote by $G_{\varepsilon}^{\pm}(z)$ the respective inverse of $T_{\varepsilon}^{\pm}(z)$. Then, there exists $C>0$, such that:
\begin{itemize}
\item For all $(\varepsilon,z) \in (0,\varepsilon_2]\times \Omega_{(\theta_0 - \delta, \theta_0 + \delta),2}^+$: $\|G_{\varepsilon}^{\pm}(z)\| \leq C \varepsilon^{-1}$.
\item For all $(\varepsilon,z) \in [0,\varepsilon_2]\times S_{(\theta_0 - \delta, \theta_0 + \delta),2}^+$: $\|G_{\varepsilon}^{\pm}(z)\| \leq C (1-|z|^2)^{-1}$.
\end{itemize}
Morever, there exists $C>0$, such that for all $(\varepsilon,z) \in (0,\varepsilon_2]\times \Omega_{(\theta_0 - \delta, \theta_0 + \delta),2}^+$ and all $\psi \in {\cal H},$
\begin{equation*}
\|G_{\varepsilon}^{\pm}(z) \psi\| \leq  C \left( \sqrt{\frac{|\bra \psi, \Re(G_{\varepsilon}^{\pm}(z)) \psi \ket |}{\varepsilon}} +\|\psi \| \right)\enspace.
\end{equation*}
\end{lem}
\noindent {\bf Proof:} Let $c=\min(c_1,c_2)$. Using Lemma \ref{lem2}, we have that for all $(\varepsilon,z) \in (0,\varepsilon_2]\times \Omega_{(\theta_0 - \delta, \theta_0 + \delta),2}^+$ (resp. $(\varepsilon,z) \in [0,\varepsilon_2]\times S_{(\theta_0 - \delta, \theta_0 + \delta),2}^+$):
\begin{eqnarray*}
(a_0 c |z|^{\pm 2} \varepsilon \mp (1-|z|^{\pm 2})) \|\psi \|^2 &\leq & \bra \psi, (c |z|^{\pm 2} \varepsilon B_1\mp (1-|z|^{\pm 2})) \psi \ket + \frac{4 c b \varepsilon}{(1-|z|)^2 +\frac{16}{\pi^2} |z| \delta^2} \|T_{\varepsilon}^{\pm}(z) \psi\|^2\\
& \leq & \pm \bra \psi, (|z|^{\pm 2} (1-e^{\mp 2 \varepsilon B_1})+(1-|z|^{\pm 2})) \psi \ket + \frac{4 c b \varepsilon}{(1-|z|)^2 +\frac{16}{\pi^2} |z| \delta^2} \|T_{\varepsilon}^{\pm}(z) \psi\|^2
\end{eqnarray*}
which leads to:
\begin{eqnarray*}
(a_0 c |z|^2 \varepsilon + (1-|z|^2)) \|\psi \|^2 &\leq & \bra \psi, (1-|z|^2(e^{-\varepsilon B(\varepsilon)})^* e^{-\varepsilon B(\varepsilon)}) \psi \ket + \frac{4 c b \varepsilon}{(1-|z|)^2 +\frac{16}{\pi^2} |z| \delta^2} \|T_{\varepsilon}^{+}(z) \psi\|^2\\
& & + |z|^2 \bra \psi, ((e^{-\varepsilon B(\varepsilon)})^* e^{-\varepsilon B(\varepsilon)}-e^{-2 \varepsilon B_1}) \psi \ket\\
(a_0 c |z|^{-2} \varepsilon - (1-|z|^{-2})) \|\psi \|^2 &\leq & \bra \psi, (|z|^{-2} e^{\varepsilon B(\varepsilon)} (e^{\varepsilon B(\varepsilon)})^*-1) \psi \ket + \frac{4 c b \varepsilon}{(1-|z|)^2 +\frac{16}{\pi^2} |z| \delta^2} \|T_{\varepsilon}^{-}(z) \psi\|^2\\
& & - |z|^{-2} \bra \psi, (e^{\varepsilon B(\varepsilon)} (e^{\varepsilon B(\varepsilon)})^*-e^{2 \varepsilon B_1}) \psi \ket
\end{eqnarray*}
Due to Lemma \ref{lem1}, there exists $\varepsilon_2 \in (0,\varepsilon_1]$ such that for all $(\varepsilon,z) \in (0,\varepsilon_2]\times \Omega_{(\theta_0 - \delta, \theta_0 + \delta),2}^+$ (resp. $(\varepsilon,z) \in [0,\varepsilon_2]\times S_{(\theta_0 - \delta, \theta_0 + \delta),2}^+$):
\begin{eqnarray*}
(\frac{a_0}{2} c |z|^2 \varepsilon + (1-|z|^2)) \|\psi \|^2 &\leq & \bra T_{\varepsilon}(z)^{+} \psi, \psi \ket + \bra \psi, \bar{z}(e^{- \varepsilon B(\varepsilon)})^* U_{\varepsilon} T_{\varepsilon}(z)^{+} \psi \ket \\
& & + \frac{4 c b \varepsilon}{(1-|z|)^2 +\frac{16}{\pi^2} |z| \delta^2} \|T_{\varepsilon}^{+}(z) \psi\|^2\\
(a_0 c |z|^{-2} \varepsilon - (1-|z|^{-2})) \|\psi \|^2 &\leq & \bra T_{\varepsilon}(z)^{-} \psi, \psi \ket + \bra \psi, z^{-1} e^{\varepsilon B(\varepsilon)} U_{\varepsilon} T_{\varepsilon}(z)^{-} \psi \ket \\
& & + \frac{4 c b \varepsilon}{(1-|z|)^2 +\frac{16}{\pi^2} |z| \delta^2} \|T_{\varepsilon}^{-}(z) \psi\|^2 \enspace.
\end{eqnarray*}
This shows that the operators  $T_{\varepsilon}^{\pm}(z)$ are injective. On the other hand, Ran $T_{\varepsilon}^{+}(z)=$ Ker $(T_{\varepsilon}^{+}(z)^*)^{\perp} ={\cal H}$ since $T_{\varepsilon}^{+}(z)^*=-\bar{z}T_{\varepsilon}^{-}(z)e^{-\varepsilon B(\varepsilon)} U$. Similarly Ran $T_{\varepsilon}^{-}(z)=$ Ker $(T_{\varepsilon}^{-}(z)^*)^{\perp} ={\cal H}$. This proves the first part of the lemma. Let $(\varepsilon,z) \in (0,\varepsilon_2]\times \Omega_{(\theta_0 - \delta, \theta_0 + \delta),2}^+$ (resp. $(\varepsilon,z) \in [0,\varepsilon_2]\times S_{(\theta_0 - \delta, \theta_0 + \delta),2}^+$), $\varphi \in {\cal H}$ and $\psi = G_{\varepsilon}^{+}(z) \varphi$. Setting $a=a_0/2$, we have by Lemma \ref{lem2}, that:
\begin{eqnarray*}
(a c |z|^{\pm 2} \varepsilon \mp (1-|z|^{\pm 2})) \|G_{\varepsilon}^{\pm}(z) \varphi \|^2 &\leq & \bra \varphi, G_{\varepsilon}^{\pm}(z) \varphi \ket - \bra \varphi, G_{\varepsilon}^{\mp}(z) \varphi \ket + \frac{4 c b \varepsilon}{(1-|z|)^2 +\frac{16}{\pi^2} |z| \delta^2} \|\varphi\|^2\\
&\leq & 2 \Re (\bra \varphi, G_{\varepsilon}^{\pm}(z) \varphi \ket) + \frac{4 c b \varepsilon}{(1-|z|)^2 +\frac{16}{\pi^2} |z| \delta^2} \|\varphi\|^2
\end{eqnarray*}
On one hand, we have that for all $(\varepsilon,z) \in (0,\varepsilon_2]\times \Omega_{(\theta_0 - \delta, \theta_0 + \delta),2}^+$ (resp. $(\varepsilon,z) \in [0,\varepsilon_2]\times S_{(\theta_0 - \delta, \theta_0 + \delta),2}^+$):
\begin{equation*}
(a c |z|^{\pm 2} \varepsilon \mp (1-|z|^{\pm 2})) \|G_{\varepsilon}^{\pm} (z) \varphi \|^2 \leq C \|\varphi\| \|G_{\varepsilon}^{\pm} (z) \varphi \| \left( 1 + \frac{\varepsilon}{(1-|z|)^2 +\frac{16}{\pi^2} |z| \delta^2}\|T_{\varepsilon}^{\pm}(z)\| \right) \enspace ,
\end{equation*}
for some $C>0$, which implies readily the first estimates (see Lemma \ref{lem1}). Since for nonnegative real numbers $x$ and $y$, $\sqrt{x+y}\leq \sqrt{x}+\sqrt{y}$, we also deduce that for all $(\varepsilon,z) \in (0,\varepsilon_2]\times \Omega_{(\theta_0 - \delta, \theta_0 + \delta),2}^+$:
\begin{equation*}
\sqrt{a c |z|^{\pm 2} \varepsilon \mp (1-|z|^{\pm 2})}\, \|G_{\varepsilon}^{\pm}(z) \varphi \| \leq C \left( |\bra \varphi, \Re (G_{\varepsilon}^{\pm}(z)) \varphi \ket|^{1/2} + \sqrt{\frac{\varepsilon}{(1-|z|)^2 +\frac{16}{\pi^2} |z| \delta^2}}\, \|\varphi\| \right) \enspace,
\end{equation*}
for some $C>0$, which implies the last estimates. \ep

So far, we have essentially dealt with the propagation properties of the operator $U$ with respect to $A$. Now, we introduce the regularity assumptions. The proof of the following lemma is straightforward:
\begin{lem}\label{lem6} Let $J\subset {\mathbb R}$ be an open bounded interval and $C$ defined by:
\begin{eqnarray*}
C: J & \rightarrow & {\cal B}({\cal H})\\
\varepsilon & \mapsto & C(\varepsilon)
\end{eqnarray*}
be a $C^1$ function with respect to the norm topology on ${\cal B}({\cal H})$. Then, the map $\varepsilon \mapsto e^{-C(\varepsilon)}$ is also norm-$C^1$ on the interval $J$. Moreover, for all $\varepsilon \in J,$
\begin{equation*}
e^{C(\varepsilon)} \partial_{\varepsilon}e^{-C(\varepsilon)}= -\sum_{p=1}^{\infty}\frac{1}{p!}\mathrm{ad}_{C(\varepsilon)}^{p-1} (\partial_{\varepsilon}{C(\varepsilon)}) \enspace .
\end{equation*}
\end{lem}

It follows that:
\begin{lem}\label{lem5} Suppose that the maps defined on $(0,\varepsilon_0]$, $\varepsilon \mapsto U_{\varepsilon}$ and $\varepsilon \mapsto B(\varepsilon)$ are $C^1$ with respect to the norm topology on ${\cal B}({\cal H})$. Then for any fixed $z \in \overline{\mathbb D} \setminus \{0\}$, the maps $\varepsilon \mapsto T_{\varepsilon}^{\pm}(z)$ are $C^1$ on $(0,\varepsilon_0]$ with respect to the norm topology on ${\cal B}({\cal H})$. Moreover, if for any $\varepsilon \in (0,\varepsilon_0]$, $U_{\varepsilon}$ and $B(\varepsilon)$ belong to $C^1(A)$ then, given $(\varepsilon,z)\in (0,\varepsilon_0]\times {\mathbb D}^*$, $T_{\varepsilon}^{\pm}(z)$ belongs to $C^1(A)$ and we have that:
\begin{eqnarray*}
\partial_{\varepsilon} T_{\varepsilon}^{+}(z) &=& -z (\partial_{\varepsilon} U_{\varepsilon}^*) e^{-\varepsilon B(\varepsilon)} -z U_{\varepsilon}^* (\partial_{\varepsilon} e^{-\varepsilon B(\varepsilon)})\\
\partial_{\varepsilon} T_{\varepsilon}^{-}(z) &=& -\bar{z}^{-1} (\partial_{\varepsilon} U_{\varepsilon}^*) e^{\varepsilon B(\varepsilon)^*} - \bar{z}^{-1} U_{\varepsilon}^* (\partial_{\varepsilon} e^{\varepsilon B(\varepsilon)^*})\\
\mathrm{ad}_A T_{\varepsilon}^{+}(z) &=& -z U_{\varepsilon}^* \left( (U_{\varepsilon} A U_{\varepsilon}^* -A)e^{-\varepsilon B(\varepsilon)} + \mathrm{ad}_A e^{-\varepsilon B(\varepsilon)} \right)\\
\mathrm{ad}_A T_{\varepsilon}^{-}(z) &=& -\bar{z}^{-1} U_{\varepsilon}^* \left( (U_{\varepsilon} A U_{\varepsilon}^* -A)e^{\varepsilon B(\varepsilon)^*} + \mathrm{ad}_A e^{\varepsilon B(\varepsilon)^*} \right) \enspace.
\end{eqnarray*}
\end{lem}
\noindent{\bf Proof:} The regularity of the maps $\varepsilon \mapsto T_{\varepsilon}^{\pm}(z)$ follows from the hypotheses and Lemma \ref{lem6} and implies naturally the formulas for $\partial_{\varepsilon} T_{\varepsilon}^{\pm}(z)$. Now, given $(\varepsilon,z)\in (0,\varepsilon_0]\times \overline{\mathbb D} \setminus \{0\}$, it follows from the hypotheses and Corollary \ref{8-2} that the operators $e^{-\varepsilon B(\varepsilon)}$ and $e^{\varepsilon B(\varepsilon)^*}$ belong to $C^1(A)$. The rest of the proof follows from Proposition \ref{3}. \ep

Note that for all $\varepsilon \in (0,\varepsilon_1],$
\begin{eqnarray*}
\partial_{\varepsilon} e^{-\varepsilon B(\varepsilon)} &=& \sum_{k=1}^{\infty} \frac{(-\varepsilon)^k}{k!} \partial_{\varepsilon} (B(\varepsilon))^k - B(\varepsilon) e^{-\varepsilon B(\varepsilon)} \\
\partial_{\varepsilon} e^{\varepsilon B(\varepsilon)^*} &=& \sum_{k=1}^{\infty} \frac{\varepsilon^k}{k!} \partial_{\varepsilon} (B(\varepsilon)^*)^k + B(\varepsilon)^* e^{\varepsilon B(\varepsilon)^*} \enspace.
\end{eqnarray*}
As a consequence, we deduce that:
\begin{lem}\label{lem7} Suppose that the maps defined on $(0,\varepsilon_0]$, $\varepsilon \mapsto U_{\varepsilon}$ and $\varepsilon \mapsto B(\varepsilon)$ are $C^1$ with respect to the norm topology on ${\cal B}({\cal H})$. Then for any fixed $z\in \Omega_{(\theta_0 - \delta, \theta_0 + \delta),2}^+$, the map $\varepsilon \mapsto G_{\varepsilon}^{\pm}(z)$ is $C^1$ on $(0,\varepsilon_2]$ with respect to the norm topology. Moreover, if for any $\varepsilon \in (0,\varepsilon_0]$, $U_{\varepsilon}$ and $B(\varepsilon)$ belong to $C^1(A)$ then, given $(\varepsilon,z) \in (0,\varepsilon_2]\times \Omega_{(\theta_0 - \delta, \theta_0 + \delta),2}^+$, $G_{\varepsilon}^{\pm}(z)$ belongs to $C^1(A)$ and we have that: $\partial_{\varepsilon} G_{\varepsilon}^{\pm}(z) = \pm \mathrm{ad}_A G_{\varepsilon}^{\pm}(z) + G_{\varepsilon}^{\pm}(z) Q^{\pm}(\varepsilon,z) G_{\varepsilon}^{\pm}(z)$ where,
\begin{eqnarray*}
Q^+(\varepsilon,z) &=& z (\partial_{\varepsilon} U_{\varepsilon}^*) e^{-\varepsilon B(\varepsilon)} + z U_{\varepsilon}^* \left( \partial_{\varepsilon} e^{-\varepsilon B(\varepsilon)} + (A-U_{\varepsilon} A U_{\varepsilon}^*) e^{-\varepsilon B(\varepsilon)} - \mathrm{ad}_A e^{-\varepsilon B(\varepsilon)} \right) \\
Q^-(\varepsilon,z) &=& \bar{z}^{-1} (\partial_{\varepsilon} U_{\varepsilon}^*) e^{\varepsilon B(\varepsilon)^*} + \bar{z}^{-1} U_{\varepsilon}^* \left( \partial_{\varepsilon} e^{\varepsilon B(\varepsilon)^*} - (A-U_{\varepsilon} A U_{\varepsilon}^*) e^{\varepsilon B(\varepsilon)^*} + \mathrm{ad}_A e^{\varepsilon B(\varepsilon)^*} \right)
\end{eqnarray*}
\end{lem}
\noindent {\bf Proof:} Since the treatment of both cases is similar, we drop the superscript $\pm$ in the proof. We observe first that given $z\in \Omega_{(\theta_0 - \delta, \theta_0 + \delta),2}^+$, for all $(\rho,\varepsilon) \in (0,\epsilon_2]^2,$
\begin{equation*}
G_{\rho}(z)-G_{\varepsilon}(z) = G_{\rho}(z) \left( T_{\varepsilon}(z)-T_{\rho}(z) \right) G_{\varepsilon}(z)
\end{equation*}
Due to Lemmatas \ref{lem4} and \ref{lem5}, $\|G_{\varepsilon}^+(z)\| \leq C\varepsilon^{-1}$ for all $\varepsilon \in (0,\varepsilon_2]$. and the map $\varepsilon \mapsto T_{\varepsilon}(z)$ is $C^1$ on $(0,\varepsilon_2]$ with respect to the norm topology. It follows that the map $\varepsilon \mapsto G_{\varepsilon}(z)$ is $C^1$ on $(0,\varepsilon_2]$ with respect to the norm topology and that:
\begin{equation*}
\partial_{\varepsilon}G_{\varepsilon}(z) = - G_{\varepsilon}(z) \left( \partial_{\varepsilon}T_{\varepsilon}(z) \right) G_{\varepsilon}(z) \enspace .
\end{equation*}
Now, fix $(\varepsilon,z) \in (0,\varepsilon_2]\times \Omega_{(\theta_0 - \delta, \theta_0 + \delta),2}^+$. It follows from Lemma \ref{lem5} and Proposition \ref{3} that $G_{\varepsilon}(z)$ belongs to $C^1(A)$ and that:
\begin{equation*}
\mathrm{ad}_A G_{\varepsilon}(z) = - G_{\varepsilon}(z) \left( \mathrm{ad}_A T_{\varepsilon}(z) \right) G_{\varepsilon}(z) \enspace .
\end{equation*}
The last part of the lemma follows by direct computation. \ep

Let us introduce more notations. Given any family of vectors $(\varphi_{\varepsilon})_{\varepsilon \in (0,\varepsilon_2]} \subset {\cal D}(A)$ such that the map $\varepsilon \mapsto \varphi_{\varepsilon}$ is $C^1$, we define the complex-valued functions $F^{\pm}$ on $(0,\varepsilon_2]\times \Omega_{(\theta_0 - \delta, \theta_0 + \delta),2}^+$ by:
\begin{equation*}
F^{\pm}(\varepsilon,z) = \bra \varphi_{\varepsilon} G_{\varepsilon}^{\pm}(z) \varphi_{\varepsilon}\ket \enspace.
\end{equation*}
\begin{lem}\label{diffin} Suppose that the maps defined on $(0,\varepsilon_0]$, $\varepsilon \mapsto U_{\varepsilon}$ and $\varepsilon \mapsto B(\varepsilon)$ are $C^1$ with respect to the norm topology on ${\cal B}({\cal H})$ and that for any $\varepsilon \in (0,\varepsilon_0]$, $U_{\varepsilon}$ and $B(\varepsilon)$ belong to $C^1(A)$. Then, for any $z\in \Omega_{(\theta_0 - \delta, \theta_0 + \delta),2}^+$, the maps $\varepsilon \mapsto F^{\pm}(\varepsilon,z)$ are of class $C^1$ and:
\begin{equation}\label{derivativeF}
\partial_{\varepsilon}F^{\pm}(\varepsilon,z)= \bra \partial_{\varepsilon}\varphi_{\varepsilon}\pm A\varphi_{\varepsilon}, G_{\varepsilon}^{\pm}(z)\varphi_{\varepsilon} \ket + \bra G_{\varepsilon}^{\pm}(z)^*\varphi_{\varepsilon}, \partial_{\varepsilon}\varphi_{\varepsilon} \mp A\varphi_{\varepsilon} \ket + \bra G_{\varepsilon}^{\pm}(z)^*\varphi_{\varepsilon}, Q^{\pm}(\varepsilon,z) G_{\varepsilon}^{\pm}(z)\varphi_{\varepsilon} \ket
\end{equation}
It follows that there exists $C>0$ such that for all $(\varepsilon,z) \in (0,\varepsilon_2]\times \Omega_{(\theta_0 - \delta, \theta_0 + \delta),2}^+,$
\begin{eqnarray}\label{diffineq}
|\partial_{\varepsilon} F^{\pm}(\varepsilon,z) | &\leq&  C\, \varepsilon q(\varepsilon) \left( \sqrt{\frac{|F^{\pm}(\varepsilon,z) |}{\varepsilon}}+\|\varphi_{\varepsilon} \| \right) \left( \sqrt{\frac{|F^{\mp}(\varepsilon,z) |}{\varepsilon}}+\|\varphi_{\varepsilon} \| \right)\nonumber \\
&+& C\, l(\varepsilon) \left( \sqrt{\frac{|F^{\pm}(\varepsilon,z) |}{\varepsilon}}+\sqrt{\frac{|F^{\mp}(\varepsilon,z) |}{\varepsilon}}+\|\varphi_{\varepsilon} \| \right)
\end{eqnarray}
with $q(\varepsilon)= \varepsilon^{-1} \max (\sup_{z \in \Omega_{(\theta_0 - \delta, \theta_0 + \delta),2}^+}\|Q^{\pm}(\varepsilon,z)\|)$ and $l(\varepsilon)=\|\partial_{\varepsilon}\varphi_{\varepsilon}\|+\|A\varphi_{\varepsilon}\|$.
\end{lem}
\noindent {\bf Proof:} Let $(\varphi_1,\varphi_2) \in {\cal D}(A)^2$. It follows from Lemma \ref{lem7} that:
\begin{eqnarray*}
\bra \varphi_1, \partial_{\varepsilon} G_{\varepsilon}^{\pm}(z) \varphi_2 \ket &=& \pm \bra A\varphi_1, G_{\varepsilon}^{\pm}(z) \varphi_2 \ket \mp \bra G_{\varepsilon}^{\pm}(z)^* \varphi_1, A\varphi_2 \ket  + \bra G_{\varepsilon}^{\pm}(z)^* \varphi_1, Q^{\pm}(\varepsilon,z) G_{\varepsilon}^{\pm}(z) \varphi_2 \ket \\
\mathrm{where} \quad G_{\varepsilon}^+(z)^* &=& - \bar{z}^{-1}U_{\varepsilon}^* e^{\varepsilon B(\varepsilon)^*} G_{\varepsilon}^-(z) \\
G_{\varepsilon}^-(z)^* &=& - z U_{\varepsilon}^* e^{\varepsilon B(\varepsilon)} G_{\varepsilon}^+(z) \enspace,
\end{eqnarray*}
for all $(\varepsilon,z) \in (0,\varepsilon_2]\times \Omega_{(\theta_0 - \delta, \theta_0 + \delta),2}^+$. This implies identity (\ref{derivativeF}). On the other hand, by Lemma \ref{lem4}, there exists $C>0$, such that for all $(\varepsilon,z) \in (0,\varepsilon_2]\times \Omega_{(\theta_0 - \delta, \theta_0 + \delta),2}^+$ and any $\psi \in {\cal H},$
\begin{equation*}
\|G_{\varepsilon}^{\pm}(z) \psi\| \leq  C \left( \sqrt{\frac{|\bra \psi, \Re (G_{\varepsilon}^{\pm}(z)) \psi \ket |}{\varepsilon}} +\|\psi \| \right)\enspace.
\end{equation*}
Inequality (\ref{diffineq}) follows. \ep

The next step consists in integrating the differential inequality of Lemma \ref{diffin}. This is done by using the following avatar of the Gronwall Lemma:
\begin{lem}\label{gronwall} Let $J=(a,b) \subset {\mathbb R}$ be an open interval and let $f$, $\varphi$ and $\psi$ be non-negative real functions on $J$ with $f$ bounded, $\varphi$ and $\psi$ in $L^1(J)$. Assume there exists $\omega \geq 0$ and $\theta \in [0,1)$ such that for all $\lambda\in J$:
\begin{equation*}
f(\lambda) \leq \omega + \int_{\lambda}^b (\varphi(\tau) f(\tau)^{\theta}+ \psi(\tau) f(\tau))\, d\tau
\end{equation*}
Then for any $\lambda\in J$, one has
\begin{equation*}
f(\lambda) \leq \left[\omega^{1-\theta}+(1-\theta) \int_{\lambda}^b \varphi(\mu) e^{(\theta-1)\int_{\mu}^b \psi(\tau)\, d\tau}\,d\mu \right]^{1/(1-\theta)} \cdot e^{\int_{\lambda}^b \psi(\tau)\, d\tau}
\end{equation*}
\end{lem}
 We refer to \cite{abmg} Appendix 7.A or \cite{hm} chapter III for a proof.

In order to apply successfully Lemma \ref{gronwall}, let us choose $\varphi \in {\cal K}$ where the interpolation space ${\cal K}:=({\cal D}(A), {\cal H})_{1/2,1}$ is continuously and densely embedded in ${\cal H}$. For such a vector $\varphi$, there exists a family of vectors $(\varphi_{\varepsilon})_{\varepsilon \in (0,\varepsilon_2]} \subset {\cal D}(A)$ such that the map $\varepsilon \mapsto \varphi_{\varepsilon}$ is $C^1$ and $\lim_{\varepsilon \rightarrow 0^+}\varphi_{\varepsilon} =\varphi$. Actually, this construction can be explicited: $\varphi_{\varepsilon}=(I+i\varepsilon A)^{-1}\varphi$ (see also \cite{abmg} Proposition 2.7.2). In this case, the function $\varepsilon \mapsto \varepsilon^{-1/2} l(\varepsilon)$ is integrable, which also implies the integrability of the function $l$. Since for all $\varepsilon \in (0, \varepsilon_2],$
\begin{equation*}
\|\varphi_{\varepsilon}\| \leq \|\varphi_{\varepsilon_2}\|+ \int_{\varepsilon_2}^{\varepsilon}\|\partial_{\tau}\varphi_{\tau}\|\, d\tau \enspace ,
\end{equation*}
the functions $\varepsilon \mapsto \|\varphi_{\varepsilon}\|$ and $\varepsilon \mapsto l(\varepsilon)\|\varphi_{\varepsilon}\|$ are also integrable.
As a consequence, we obtain:
\begin{lem}\label{gronwall-1} Let $\varphi \in {\cal K}$ and fix a family of vectors $(\varphi_{\varepsilon})_{\varepsilon \in (0,\varepsilon_2]} \subset {\cal D}(A)$ such that the map $\varepsilon \mapsto \varphi_{\varepsilon}$ is $C^1$ and $\lim_{\varepsilon \rightarrow 0^+}\varphi_{\varepsilon} =\varphi$. Suppose that the maps defined on $(0,\varepsilon_0]$, $\varepsilon \mapsto U_{\varepsilon}$ and $\varepsilon \mapsto B(\varepsilon)$ are $C^1$ with respect to the norm topology on ${\cal B}({\cal H})$ and that for any $\varepsilon \in (0,\varepsilon_0]$, $U_{\varepsilon}$ and $B(\varepsilon)$ belong to $C^1(A)$. If
\begin{equation*}
\int_0^{\varepsilon_0} q(\varepsilon)\, d\varepsilon < \infty \enspace ,
\end{equation*}
then there exist $C>0$ and $H\in L^1((0,\varepsilon_2])$ such that for all $(\varepsilon,z) \in (0,\varepsilon_2]\times \Omega_{(\theta_0 - \delta, \theta_0 + \delta),2}^+,$
\begin{eqnarray*}
|F^{\pm}(\varepsilon,z)| &<& C \\
|\partial_{\varepsilon}F^{\pm}(\varepsilon,z)| &\leq & H(\varepsilon) \enspace.
\end{eqnarray*}
\end{lem}
\noindent {\bf Proof:} The reader will observe first that the integrability of the function $q$ implies the integrability of the function $\varepsilon \mapsto \varepsilon q(\varepsilon)\|\varphi_{\varepsilon} \|^2$. Define, the auxiliary functions $K$ and $L$ by
\begin{eqnarray*}
K(\varepsilon,z) &=& | F^+(\varepsilon,z)| +| F^-(\varepsilon,z)| \\
L(\varepsilon) &=& \sup_{z\in \Omega_{(\theta_0 - \delta, \theta_0 + \delta),2}^+} K(\varepsilon,z)
\end{eqnarray*}
Up some adjustment of the constants, we have that for all $(\varepsilon,z) \in (0,\varepsilon_2]\times \Omega_{(\theta_0 - \delta, \theta_0 + \delta),2}^+,$
\begin{eqnarray*}
\left| K(\varepsilon_2,z) - K(\varepsilon,z) \right| &= & \left| | F^+(\varepsilon_2,z)| - | F^+(\varepsilon,z)|+  | F^-(\varepsilon_2,z)| -| F^-(\varepsilon,z)| \right| \\
&\leq & | F^+(\varepsilon_2,z) - F^+(\varepsilon,z)|+  | F^-(\varepsilon_2,z) - F^-(\varepsilon,z)| \\
& \leq & \int_{\varepsilon}^{\varepsilon_2} | \partial_{\rho} F^+(\rho,z) |+| \partial_{\rho} F^-(\rho,z) | \, d\rho \\
& \leq & C \int_{\varepsilon}^{\varepsilon_2} (q(\rho)  K(\rho,z) +l(\rho) \rho^{-1/2} K(\rho,z)^{1/2} +\rho q(\rho) \|f_{\rho}\|^2 +l(\rho) \|f_{\rho}\|) \, d\rho 
\end{eqnarray*}
using Lemma \ref{diffin} and the fact that: $| F^{\pm}(\varepsilon,z)| \leq K(\varepsilon,z)$. It follows from Lemma \ref{lem4} that for all $(\varepsilon,z) \in (0,\varepsilon_2]\times \Omega_{(\theta_0 - \delta, \theta_0 + \delta),2}^+,$
\begin{eqnarray*}
K(\varepsilon,z) &\leq & K(\varepsilon_2,z) + C \int_{\varepsilon}^{\varepsilon_2} (q(\rho)  K(\rho,z) +l(\rho)\rho^{-1/2} K(\rho,z)^{1/2} +\rho q(\rho) \|f_{\rho}\|^2 +l(\rho) \|f_{\rho}\|) \, d\rho \\
&\leq & C\left(\varepsilon_2^{-1} +\int_{\varepsilon}^{\varepsilon_2} (q(\rho)  K(\rho,z) +l(\rho)\rho^{-1/2} K(\rho,z)^{1/2} +\rho q(\rho) \|f_{\rho}\|^2 +l(\rho) \|f_{\rho}\|) \, d\rho \right) \\
L(\varepsilon) &\leq & C\left( \varepsilon_2^{-1} + \int_{\varepsilon}^{\varepsilon_2} (q(\rho)  L(\rho) +l(\rho)\rho^{-1/2} L(\rho)^{1/2} +\rho q(\rho) \|f_{\rho}\|^2 +l(\rho) \|f_{\rho}\|) \, d\rho \right) \enspace .
\end{eqnarray*}
The first estimate follows from Lemma \ref{gronwall}. The second part is obtained, plugging the first estimate in the differential inequality (\ref{diffineq}). \ep

Let us explicit the implication of Lemma \ref{gronwall-1}:
\begin{cor}\label{cor1} Let $(\varphi,\psi)\in {\cal K}^2$ and denote $\Theta_0:=(\theta_0-\delta,\theta_0+\delta)$. Under the hypotheses of Lemma \ref{gronwall-1}, the holomorphic functions defined on ${\mathbb D}$ by $z\mapsto \bra \varphi, (1-zU^*)^{-1}\psi\ket$ and $z\mapsto \bra \varphi, (1-\bar{z}^{-1}U^*)^{-1}\psi\ket$ extend continuously to ${\mathbb D}\cup e^{i\Theta_0}$.
\end{cor}
\noindent{\bf Proof:} Due to the polarization identity, it is enough to prove the result for the case $\varphi=\psi$. It follows from Lemma \ref{gronwall-1} that the limits $\lim_{\varepsilon \rightarrow 0} \bra \varphi_{\varepsilon}, G_{\varepsilon}^{\pm}(z) \varphi_{\varepsilon}\ket$ exist and satisfy:
\begin{equation*}
\lim_{\varepsilon \rightarrow 0} \bra \varphi_{\varepsilon}, G_{\varepsilon}^{\pm}(z) \varphi_{\varepsilon}\ket = F^{\pm}(\varepsilon_2,z) - \int_0^{\epsilon_2} \partial_{\rho}F^{\pm}(\rho,z)\, d\rho \enspace ,
\end{equation*}
for any $z\in \Omega_{(\theta_0 - \delta, \theta_0 + \delta),2}^+$. On the other hand, for any $z\in {\mathbb D}$, $\lim_{\varepsilon \rightarrow 0^+} T_{\varepsilon}^+(z)=(1-zU^*)$ and $\lim_{\varepsilon \rightarrow 0^+} T_{\varepsilon}^-(z)=(1-\bar{z}^{-1}U^*)$ in the norm topology of ${\cal B}({\cal H})$. It results from the resolvent identity that $\lim_{\varepsilon \rightarrow 0^+} G_{\varepsilon}^+(z)=(1-zU^*)^{-1}$ and $\lim_{\varepsilon \rightarrow 0^+} G_{\varepsilon}^-(z)=(1-\bar{z}^{-1}U^*)^{-1}$ w.r.t the same topology. For any $(\varepsilon,z)\in (0,\varepsilon_2]\times S_{(\theta_0 - \delta, \theta_0 + \delta),2}^+,$
\begin{eqnarray*}
| F^+(\varepsilon, z) - \bra \varphi, (1-zU^*)^{-1}\varphi\ket | &\leq & \|\varphi-\varphi_{\varepsilon}\| \left(\|G^+_{\varepsilon}(z)\|+\|G^+_0(z)\|\right) \|\varphi_{\varepsilon}\| +\|G^+_{\varepsilon}(z)-G^+_0(z)\| \|\varphi_{\varepsilon}\|^2\\
| F^-(\varepsilon, z) - \bra \varphi, (1-\bar{z}^{-1}U^*)^{-1}\varphi\ket | &\leq & \|\varphi-\varphi_{\varepsilon}\| \left(\|G^-_{\varepsilon}(z)\|+\|G^-_0(z)\|\right) \|\varphi_{\varepsilon}\| +\|G^-_{\varepsilon}(z)-G^-_0(z)\| \|\varphi_{\varepsilon}\|^2 \enspace,
\end{eqnarray*}
with $\|G^{\pm}_{\varepsilon}(z)\| \leq C (1-|z|^2)^{-1}$ for some $C>0$. This implies that for any $z\in S_{(\theta_0 - \delta, \theta_0 + \delta),2}^+,$
\begin{eqnarray*}
\bra \varphi, (1-zU^*)^{-1}\varphi\ket &=& \lim_{\varepsilon \rightarrow 0} \bra \varphi_{\varepsilon}, G_{\varepsilon}^+(z) \varphi_{\varepsilon}\ket \\
\bra \varphi, (1-\bar{z}^{-1}U^*)^{-1}\varphi\ket &=& \lim_{\varepsilon \rightarrow 0} \bra \varphi_{\varepsilon}, G_{\varepsilon}^-(z) \varphi_{\varepsilon}\ket \enspace .
\end{eqnarray*}
On the other hand, for any $\varepsilon \in (0,\varepsilon_2]$, the functions $z\mapsto G_{\varepsilon}^{\pm}(z)$ and $z\mapsto \partial_{\varepsilon} F^{\pm}(\varepsilon,z)$ are respectively norm-continuous and continuous on $\Omega_{(\theta_0 - \delta, \theta_0 + \delta),2}^+$ (see identity (\ref{derivativeF})). The desired result follows then from Lemma \ref{gronwall-1} and the dominated convergence Theorem.  \ep
\\

It remains to relate the integrability of the function $q$ to the hypothesis $U\in {\cal C}^{1,1}(A)$. The next result is the following step in this direction:
\begin{lem}\label{qepsilon} Assume there exist two families of respectively unitary and bounded operators $(B(\varepsilon))_{\varepsilon \in (0,\varepsilon_0]}$ and $(U_{\varepsilon})_{\varepsilon \in (0,\varepsilon_0]}$ on ${\cal B}({\cal H})$ such that:
\begin{itemize}
\item there exists $C>0$, such that for all $\varepsilon \in (0,\varepsilon_0]$, $\varepsilon \, \|B(\varepsilon)\| + \|U_{\varepsilon}-U\| \leq C \varepsilon$,
\item $\lim_{\varepsilon \rightarrow 0}\|B(\varepsilon)-(A-UAU^*)\| =0$,
\item the maps $\varepsilon \mapsto U_{\varepsilon}$ and $\varepsilon \mapsto B(\varepsilon)$ are $C^1$ with respect to the norm topology on ${\cal B}({\cal H})$
\item for any $\varepsilon \in (0,\varepsilon_0]$, $U_{\varepsilon}$ and $B(\varepsilon)$ belong to $C^1(A)$,
\item the map $\varepsilon \mapsto \varepsilon^{-1} \|\partial_{\varepsilon} U_{\varepsilon}\|+\|\partial_{\varepsilon} B(\varepsilon)\|+\|\mathrm{ad}_A B(\varepsilon)\| + \varepsilon^{-1} \|B(\varepsilon)-(A-U_{\varepsilon}AU_{\varepsilon}^*)\|$ belongs to $L^1(0,\varepsilon_0)$.
\end{itemize}
Then, the function $q$ defined by: $q(\varepsilon) = \varepsilon^{-1} \max (\sup_{z \in \Omega_{{\mathbb T},2}^+}\|Q^{\pm}(\varepsilon,z)\|)$ for any $\varepsilon \in (0,\varepsilon_0]$ (with $Q^{\pm}(\varepsilon,z)$ defined in Lemma \ref{lem7}) belongs to $L^1(0,\varepsilon_0)$.
\end{lem}
\noindent{\bf Proof:} The conclusion follows from the definition of $Q^{\pm}(\varepsilon,z)$, Lemma \ref{lem7}, once noted that there exists $C>0$ such that for all $\varepsilon \in (0,\varepsilon_0]$: $\|e^{-\varepsilon B(\varepsilon)} \|\leq C$, $\|e^{\varepsilon B(\varepsilon)^*} \|\leq C$, $\| (\partial_{\varepsilon} U_{\varepsilon})^* \| = \| \partial_{\varepsilon} U_{\varepsilon} \|$,
\begin{eqnarray*}
\| \sum_{k=1}^{\infty} \frac{(-\varepsilon)^k}{k!} \partial_{\varepsilon} (B(\varepsilon))^k \| &\leq & \sum_{k=1}^{\infty} \frac{\varepsilon^k}{k!} \| \partial_{\varepsilon} (B(\varepsilon))^k \| \leq C \varepsilon \| \partial_{\varepsilon} B(\varepsilon) \| \\
\|\mathrm{ad}_A e^{-\varepsilon B(\varepsilon)}\| &\leq & C \varepsilon \|\mathrm{ad}_A B(\varepsilon)\| \\
\|\mathrm{ad}_A e^{\varepsilon B(\varepsilon)^*}\| &\leq & C \varepsilon \|\mathrm{ad}_A B(\varepsilon)\|\enspace .
\end{eqnarray*}
\ep

It turns out that the hypotheses of Lemma \ref{qepsilon} are equivalent to the hypothesis $U\in {\cal C}^{1,1}(A)$ (see Proposition \ref{georgescu} and Corollary \ref{u2}).

\subsection{Proof of Theorem \ref{thm0}}

Recall that $U$ is strictly propagating with respect to the observable $A$ on the arc $e^{i\Theta}$, $\Theta \in {\cal B}({\mathbb T})$ if there exist $c >0$ such that: $E_{\Theta} (A-UAU^*) E_{\Theta} \geq c E_{\Theta}$.

In view of Corollary \ref{virial-2}, we know that $U$ has at most a finite number of eigenvalues in $\Theta$. These eigenvalues have finite multiplicity. We also know that for any $\theta \in \Theta\setminus \sigma_{pp}(U)$, there exists $\delta_{\theta}>0$ such that $U$ is strictly propagating with respect to $A$ on $(\theta-2\delta_{\theta},\theta+2\delta_{\theta})$. Given any compact subset $K\subset \Theta\setminus \sigma_{pp}(U)$, the collection $((\theta-\delta_{\theta},\theta+\delta_{\theta}))_{\theta \in \Theta\setminus \sigma_{pp}(U)}$ induces an open covering of $K$, from which we can extract a finite open covering. Due to Lemma \ref{qepsilon}, Proposition \ref{georgescu} and Corollary \ref{u2}, Corollary \ref{cor1} applies on each of these intervals, which proves statement (ii) of Theorem \ref{thm0}. The last part follows from Proposition \ref{usmoothext}.

If in addition, $U$ is strictly propagating with respect to $A$ on $\Theta$, it clearly follows from Theorems \ref{virial} and \ref{thm2} that $U$ has even no eigenvalues in $\Theta$.

\section{Regularity Classes for Bounded operators}

This section gathers some elementary properties of the regularity classes $C^k(A)$ (sometimes denoted $C^k(A,{\cal H})$ or $C^k(A,{\cal H},{\cal H})$) introduced in Section 2. For more details see \cite{abmg} Chapter 5. From now, $A$ is a fixed self-adjoint operator, densely defined on a fixed Hilbert space ${\cal H}$, with domain ${\cal D}(A)$. 

\subsection{Basics}

The regularity of a bounded operator defined on ${\cal H}$ w.r.t $A$ is associated to the algebra of derivation on ${\cal B}({\cal H})$ defined by the operation $\mathrm{ad}_A$. From a theoretical point of view, it is often more convenient to reformulate this concept of derivation in terms of the regularity of the strongly continuous function:
\begin{eqnarray*}
{\cal W}_B: {\mathbb R} & \rightarrow & {\cal B}({\cal H})\\
t & \mapsto & e^{iAt}Be^{-iAt} \enspace .
\end{eqnarray*}
Most of the properties derived below can be deduced easily once established the following equivalence:
\begin{prop}\label{equiv0} Let $k\in {\mathbb N}$. The following assertions are equivalent:
\begin{itemize}
\item $B\in C^k(A)$
\item The map ${\cal W}_B$ is $C^k$ with respect to the strong topology on ${\cal B}({\cal H})$.
\item The map ${\cal W}_B$ is $C^k$ with respect to the weak topology on ${\cal B}({\cal H})$.
\end{itemize}
Moreover, ${\cal W}_B^{(k)} (0) = i^k \mathrm{ad}_A^k B$.
\end{prop}
For a proof, see \cite{abmg} Lemma 6.2.9, Theorem 6.2.10 in association with Lemma 6.2.1 and Definition 6.2.2. For all nonnegative integral number $k$, $C^{k+1}(A) \subset C^k(A)$.

\begin{prop}\label{2} If $B\in C^1(A)$, then $B({\cal D}(A)) \subset {\cal D}(A)$.
\end{prop}
\noindent{\bf Proof:} If $B\in C^1(A)$, then there exists $C>0$ such that for all $(\varphi,\psi) \in {\cal D}(A)\times {\cal D}(A),$
\begin{equation*}
\left| \bra A\varphi,B\psi\ket -  \bra \varphi,BA\psi\ket \right| \leq C \|\varphi\|\|\psi\| \enspace,
\end{equation*}
which implies that: $| \bra A\varphi,B\psi\ket | \leq C \|\varphi\|\left(\|\psi\|+\|A\psi\|\right)$. Hence, for all $\psi \in {\cal D}(A)$, $B\psi \in {\cal D}(A^*)={\cal D}(A)$ since $A$ is self-adjoint. \ep

For any nonnegative integral number $k$, $C^k(A)$ is clearly a vector subspace of ${\cal B}({\cal H})$. These classes also share the following algebraic properties:
\begin{prop}\label{3} Let $k\in {\mathbb N}$ and $(B,C) \in C^k(A)\times C^k(A)$. then,
\begin{itemize}
\item $B^* \in C^k (A)$ and for all $j\in \{0,\ldots,k\}$, $\mathrm{ad}_A^j B^* = (-1)^j (\mathrm{ad}_A^j B)^*$
\item $BC \in C^k(A)$ and for all $j\in \{1,\ldots,k\},$
\begin{equation*}
\mathrm{ad}_A^j BC = \sum_{l_1+l_2=j} \frac{j!}{l_1! l_2!} \mathrm{ad}_A^{l_1} B \mathrm{ad}_A^{l_2} C \enspace.
\end{equation*}
In particular, $\mathrm{ad}_A BC = (\mathrm{ad}_A B) C + B (\mathrm{ad}_A C)$
\item for all $j\in \{0,\ldots, k\}$, $\mathrm{ad}_A^j B \in C^{k-j}(A)$.
\item If $B$ is invertible (i.e $B^{-1}\in {\cal B}({\cal H})$) and $B\in C^1(A)$, then $B^{-1} \in C^1(A)$: $\mathrm{ad}_A B^{-1} = - B^{-1} (\mathrm{ad}_A B) B^{-1}$.
\end{itemize}
\end{prop}
See \cite{abmg} Propositions 5.1.2, 5.1.5, 5.1.6, 5.1.7 for a proof. Combining the last sentences of Proposition \ref{3}, we deduce that if an invertible bounded operator $B$ belongs to $C^k(A)$, then its inverse $B^{-1}$ also belongs to $C^k(A)$.

As mentioned in Section 2, if a unitary operator $U$ belongs to $C^1(A)$, then $U({\cal D}(A)) ={\cal D}(A)=U^*({\cal D}(A))$. We prove Lemma \ref{equiv1}.

\noindent{\bf Proof of Lemma \ref{equiv1}:} We prove the equivalence between statements (a) and (c). The proof of the remaining parts is similar. Assume first that $U\in C^1(A)$. As observed, $U({\cal D}(A))={\cal D}(A)$, which allows us to choose ${\cal S}_1={\cal D}(A)$. Therefore, for all $(\varphi,\psi) \in {\cal D}(A)^2,$
\begin{eqnarray*}
\left| \bra U\varphi,AU\psi\ket -  \bra \varphi,A\psi\ket \right| &=& \left| \bra U\varphi,AU\psi\ket -  \bra U \varphi,UA\psi\ket \right| = |\bra U\varphi, (\mathrm{ad}_A U) \psi \ket |\\
&\leq & \|\mathrm{ad}_A (U)\| \|\psi\| \|\varphi\| \enspace.
\end{eqnarray*}
This implies (c). Now, assume (c) holds and choose ${\cal S}={\cal S}_1$. Since ${\cal S} \subset {\cal D}(A)$ and $U^*{\cal S}={\cal S}=U{\cal S}$, we have that for all $(\varphi,\psi) \in {\cal S}^2,$
\begin{eqnarray*}
\left| \bra A\varphi,U\psi\ket -  \bra \varphi,UA\psi\ket \right| &=& \left| \bra \varphi,AU\psi\ket -  \bra \varphi,UA\psi\ket \right| = |\bra U (U^*\varphi), AU\psi \ket - \bra U^*\varphi,A \psi \ket |\\
&\leq & \| U^{*}AU-A \| \|\psi\| \|\varphi\| \enspace.
\end{eqnarray*}
which implies (a), in view of Definition \ref{c1}. \ep

As a consequence of Proposition \ref{3} and Lemma \ref{equiv1} gives:
\begin{lem}\label{equiv2} Let $k\in {\mathbb N}$ and $U$ be a unitary operator on ${\cal H}$. Then, the three following assertions are equivalent:
\begin{itemize}
\item[(a)] $U\in C^k(A)$.
\item[(b)] $U^*\in C^k(A)$.
\item[(c)] $U$ satisfies item (c) of Lemma \ref{equiv1} and $(U^* AU-A) \in C^{k-1}(A)$.
\item[(d)] $U$ satisfies item (d) of Lemma \ref{equiv1} and $(A-UAU^*) \in C^{k-1}(A)$.
\end{itemize}
\end{lem}
\noindent {\bf Proof:} If $k=0$, this is Lemma \ref{equiv1}. If $U\in C^k(A)$ then $U^*\in C^k(A)$ and $\mathrm{ad}_A U \in C^{k-1}(A)$ (see Proposition \ref{3}). Since  $U^* AU-A=U^* (\mathrm{ad}_A U)$, then $(U^* AU-A) \in C^{k-1}(A)$. The converse implication is proven by induction on $j$, $j\in \{1,\ldots,k\}$. Assume that $U\in C^j(A)$ for some $j\leq k-1$. Since $(U^* AU-A) \in C^{k-1}(A) \subset C^j(A)$ and $\mathrm{ad}_A U = U(U^* AU-A)$ it follows from Proposition \ref{3} that: $\mathrm{ad}_A U \in C^j(A)$, which means that $U\in C^{j+1}(A)$ in view of the induction hypothesis. This proves the equivalence between (a) and (c). The remaining parts can be justified similarly. \ep

\begin{prop}\label{6-1} Let $k\in {\mathbb N}$, $U$ and $B$ be respectively a unitary operator and a bounded operator defined on ${\cal H}$. Then, $U^*BU \in C^k(A)$ if and only if $B \in C^k(UAU^*)$. Moreover, for all $j\in \{0,\ldots,k\}$, $\mathrm{ad}_A^j (U^*BU) =U^* (\mathrm{ad}_{UAU^*}^j B)U$.
\end{prop}
\noindent{\bf Proof:} This follows from Proposition \ref{equiv0} and the fact that for all $t\in {\mathbb R},$
\begin{equation*}
e^{iAt}(U^*BU)e^{-iAt} = U^*\left(e^{iUAU^*t}Be^{-iUAU^*t}\right) U \enspace .
\end{equation*}
\ep

As mentioned in Section 2, it also possible to define intermediate scale of regularity. This leads us to the definition of the classes ${\cal C}^{s,p}(A)$. For more details, see \cite{abmg} paragraph 5.2 or \cite{sah}:
\begin{defin}\label{csp} Let $B\in {\cal B}({\cal H})$, $s>0$ and $p\in [1,\infty)$. Then, $B$ is of class ${\cal C}^{s,p}(A)$ if there exists $l>s$ such that:
\begin{equation}\label{int}
\int_{-1}^1 \| \sum_{m=0}^l (-1)^m \binom{l}{m}e^{imA\tau}Be^{-imA\tau} \|^p \, \frac{d\tau}{|\tau|^{sp+1}}  < \infty \enspace.
\end{equation}
\end{defin}
The classes ${\cal C}^{s,p}(A)$ are clearly vector subspaces of ${\cal B}({\cal H})$ which share the following algebraic properties:
\begin{prop}\label{4} Let $(s,p)\in {\mathbb N}$ and $(B,C) \in {\cal C}^{s,p}(A)\times {\cal C}^{s,p}(A)$. Then,
\begin{itemize}
\item $B^* \in {\cal C}^{s,p}(A)$
\item $BC \in {\cal C}^{s,p}(A)$
\item If $s=k+\sigma$ with $k\in {\mathbb N}$, $\sigma\in (0,1]$ and $j\in \{0,\ldots, k\}$, $\mathrm{ad}_A^j B \in {\cal C}^{s-j,p}(A)$.
\item If $B$ is invertible (i.e $B^{-1}\in {\cal B}({\cal H})$), then $B^{-1} \in {\cal C}^{s,p}(A)$.
\end{itemize}
\end{prop}
See \cite{abmg} Propositions 5.2.2, 5.2.3, 5.2.4 for a proof.

For the relationships between these regularity classes with the self-adjoint functional calculus, we refer the reader to \cite{abmg} Theorem 6.2.5 and Corollary 6.2.6 or to \cite{gj} and the Helffer-Sjostrand formula. related aspects are considered in the next paragraph.

\subsection{Functional Calculus}

We start by the following fundamental lemma:
\begin{lem}\label{H0-1} Let $k\in {\mathbb N}$, $(B_n)_{n\in {\mathbb N}}\subset C^1(A)$. Assume that the sequences $(B_n)$ and $(\mathrm{ad}_A B_n)$ are strongly convergent to some (bounded) operators $C_0$ and $C_1$ respectively. Then, $C_0 =s-\lim B_n$ belongs to $C^1(A)$ and
\begin{equation*}
i \mathrm{ad}_A C_0 =C_1 \enspace .
\end{equation*}
\end{lem}
\noindent{\bf Proof.} Let $\varphi \in {\cal H}$, $\|\varphi \|=1$ and ${\cal W}_{B_n}\varphi: {\mathbb R}\rightarrow {\cal H}$ defined by ${\cal W}_{B_n}(t)\varphi= e^{iAt}B_n e^{-iAt}\varphi$. By hypothesis, the functions ${\cal W}_{B_n}\varphi$ are $C^1$ and
\begin{equation*}
{\cal W}_{B_n}'(t)\varphi = i e^{iAt}(\mathrm{ad}_A B_n) e^{-iAt}\varphi
\end{equation*}
Moreover, the sequences $({\cal W}_{B_n}\varphi)$ and $({\cal W}_{B_n}'\varphi)$ converge pointwise respectively to some functions ${\cal W}_{C_0}\varphi$ and ${\cal W}_{C_1}\varphi$. We show that this convergence is actually uniform on each compact subset of ${\mathbb R}$. So, fix $a<b$ and $\epsilon >0$. The set $\{e^{-iAt}\varphi; t\in [a,b]\}$ is a compact subset of the unit sphere of ${\cal H}$. Therefore, there exists a finite rank operator $P$ such that: for all $t\in [a,b]$, $\|(I-P)e^{-iAt}\varphi\| \leq \epsilon$. On the other and, for all $t\in [a,b],$
\begin{eqnarray*}
{\cal W}_{B_n}(t)\varphi-{\cal W}_{C_0}(t)\varphi &=& e^{iAt}(B_n - C_0) P e^{-iAt}\varphi + e^{iAt}(B_n - C_0) (I-P) e^{-iAt}\varphi\\
{\cal W}_{B_n}'(t)\varphi-{\cal W}_{C_1}\varphi(t) &=& e^{iAt}(i\mathrm{ad}_A B_n - C_1) P e^{-iAt}\varphi + e^{iAt}(i\mathrm{ad}_A B_n - C_1) (I-P) e^{-iAt}\varphi \enspace.
\end{eqnarray*}
Since $P$ is a compact operator, the sequences $((\mathrm{ad}_A^j B_n) P)$ converge in norm to the bounded operators $C_jP$, i.e: there exists $N\in {\mathbb N}$ such that for all $n\geq N$ and all $t\in [a,b],$
\begin{eqnarray*}
\| e^{iAt}(B_n - C_0) P e^{-iAt}\varphi\| &\leq & \|(B_n - C_0) P\|\|\varphi\| \leq \epsilon \\
\| e^{iAt}(i\mathrm{ad}_A B_n - C_1) P e^{-iAt}\varphi\| &\leq & \|(\mathrm{ad}_A B_n - C_1) P\| \|\varphi\|\leq \epsilon \enspace.
\end{eqnarray*}
Aside from the fact that the functions ${\cal W}_{C_0}\varphi$ and ${\cal W}_{C_1}\varphi$ are continuous on $[a,b]$ we also have that ${\cal W}_{C_0}\varphi$ is of class $C^1$ on $(a,b)$ and that ${\cal W}_{C_0}'\varphi={\cal W}_{C_1}\varphi$. Since the choice of $a$, $b$ and $\varphi$ was arbitrary, the lemma is proved. \ep

It follows by induction that:
\begin{lem}\label{H0-2} Let $k\in {\mathbb N}$, $(B_n)_{n\in {\mathbb N}}\subset C^k(A)$. Assume that for all $j\in \{0,\ldots,k\}$, the sequences $(\mathrm{ad}_A^j B_n)$ converge strongly to some (bounded) operator $C_j$ respectively. Then, $C_0 =s-\lim_{n\rightarrow \infty} B_n$ belongs to $C^k(A)$ and
\begin{equation*}
\mathrm{ad}_A^j C_0 =C_j \enspace .
\end{equation*}
\end{lem}

The next result can be deduced now by considering Riemann sums:
\begin{prop}\label{8} Let $k\in {\mathbb N}$ and $(B(t))_{t\in [a,b]}$ ($-\infty<a<b<\infty$) a family of bounded operators such that:
\begin{itemize}
\item For all $t\in [a,b]$, $B(t) \in C^k(A)$.
\item For all $j\in \{0,\ldots,k\}$, the maps defined on $[a,b]$ by $t\mapsto \mathrm{ad}_A^j B(t)$ are strongly continuous.
\end{itemize}
Then, the strongly convergent integral $\int_a^b B(t)\, dt$ belongs to $C^k(A)$ and for all $j\in \{0,\ldots,k\}$
\begin{equation*}
\mathrm{ad}_A^j (\int_a^b B(t)\, dt) = \int_a^b \mathrm{ad}_A^j (B(t))\, dt \enspace.
\end{equation*}
\end{prop}


Before going further, let us recall that if $\Phi$ is a complex-valued holomorphic function on some open simply connected domain ${\cal D} \subset {\mathbb C}$, $\sigma(B) \subset {\cal D}$ and $\Gamma$ is a rectifiable Jordan curve in ${\cal D}$, which contains $\sigma(B)$ in its interior then,
\begin{equation}\label{riesz}
\Phi(B) =\frac{1}{2i\pi} \int_{\Gamma} \Phi(z) (z-B)^{-1}\, dz
\end{equation}
where the RHS is a norm convergent integral \cite{du}, \cite{Kato1}, \cite{D}. It follows that:
\begin{cor}\label{8-1} Let $k\in {\mathbb N}$. Let $\Phi$ be a complex-valued holomorphic function on some open simply connected domain ${\cal D} \subset {\mathbb C}$. If $B\in C^k(A)$ and $\sigma(B) \subset {\cal D}$, then $\Phi(B) \in C^k(A)$ and
for all $j\in \{0,\ldots,k\}$
\begin{equation*}
\mathrm{ad}_A^j (\Phi(B)) = \frac{1}{2i \pi}\int_{\Gamma} \Phi(z) \mathrm{ad}_A^j (z-B)^{-1}\, dz \enspace.
\end{equation*}
\end{cor}
\noindent {\bf Proof:} The proof follows from Proposition \ref{8} and formula (\ref{riesz}) by considering a suitable parametrization of the curve $\Gamma$. \ep

\begin{cor}\label{8-2}{\bf Baker-Campbell-Hausdorff} Let $k\in {\mathbb N}$. If $B \in C^k(A)$, then for any $\mu \in {\mathbb C}$, $e^{\mu B} \in C^k(A)$. Moreover,
\begin{equation*}
e^{-B} A e^{B}-A=\sum_{k=1}^{\infty}  \frac{(-1)^{k-1}}{k!}\, \mathrm{ad}_B^{k-1}(\mathrm{ad}_A B) \enspace.
\end{equation*}
\begin{eqnarray*}
\| \mathrm{ad}_A e^{iB} \| &\leq & e^{\|B\|} \| \mathrm{ad}_A B\| \\
\| e^{-B} A e^{B}-A \| &\leq & e^{\|B\|}\| \mathrm{ad}_A B \|\enspace .
\end{eqnarray*}
\end{cor}

Another useful consequence of Proposition \ref{8} is:
\begin{cor}\label{dyson}{\bf Dyson expansions.} Let $k\in {\mathbb N}$, $T>0$ and $(V(t,s))_{(s,t)\in [0,T]^2}$ be a strongly continuous family of bounded operators such that for all $(s,t)\in [0,T]^2$, $V(t,s)\in C^k(A)$. Then, the for any $(s,t)\in [0,T]^2$, the operator $\Omega (t,s)$ defined by the norm convergent series :
\begin{equation}\label{dyson2}
\Omega(t,s) = I + \sum_{j=1}^{\infty} (-i)^j \int_s^t \ldots \int_s^{\tau_{j-1}} V(\tau_1,s) \dots V(\tau_j,s) \, d\tau_j \ldots d\tau_1 \enspace,
\end{equation}
belongs to $C^k(A)$.
\end{cor}
We refer to \cite{rs} Theorem X.69 for more details on the Dyson expansions.

The regularity of a unitary operator has some specific features:
\begin{lem} \label{9} Let $k\in {\mathbb N}$. Assume that $U \in C^k(A)$. Then, for all $m\in {\mathbb Z}$, $U^m \in C^k(A)$. Moreover, there exists $C>0$ such that for all $j\in \{0,\ldots,k\}$ and all $|m| \geq j,$
\begin{equation*}
\| \mathrm{ad}_A^j U^m\| \leq C^j |m|^j \enspace.
\end{equation*}
Actually, $C=\|U\|_{C^k} :=\sqrt{\sum_{j=0}^{k} \|\mathrm{ad}_A^j U \|^2}\geq \|\mathrm{ad}_A^0 U \|=\|U \|=1$.
\end{lem}
\noindent{\bf Proof:} The first part is a consequence of Proposition \ref{3}. To prove the second part, it is enough to consider the case of positive $m$ (see again Proposition \ref{3}). We have that:
\begin{equation*}
\mathrm{ad}_A^j U^m = \sum_{(j_1,\ldots,j_m)\in \{0,\ldots,j\}^m, j_1+\ldots+j_m=j} \binom{j}{j_1 \ldots j_m} (\mathrm{ad}_A^{j_1} U)\ldots (\mathrm{ad}_A^{j_m} U) \enspace.
\end{equation*}
If $m\geq j$, then for each term involved in the sum on the RHS, at least $(m-j)$ coefficients $j_m$'s are zero. Since $\|\mathrm{ad}_A^0 U \|=\|U \|=1$, each of these can be estimated by $\|U\|_{C^k}^l$. Since the sum involves $m^j$ terms, the estimate follows. \ep

\begin{prop}\label{10} Let $k\in {\mathbb N}$. Assume that $U \in C^k(A)$ and $\Phi \in C^{k+2}({\mathbb T})$. Then  $\Phi(U) \in C^k(A)$ and for all $j\in \{0,\ldots,k\},$
\begin{equation*}
\mathrm{ad}_A^j \Phi(U) = \sum_{m\in {\mathbb Z}} \hat{\phi}_m \mathrm{ad}_A^j U^m \enspace.
\end{equation*}
\end{prop}
\noindent{\bf Proof:} We have that: $\Phi(U) =\sum_{m\in {\mathbb Z}} \hat{\phi}_m U^m$ where the series on the RHS is norm convergent. Actually, if $\Phi \in C^{k+2}({\mathbb T})$, then $(m^k \hat{\phi}_m)_{m\in {\mathbb Z}}\in l^1({\mathbb Z})$. It follows from Lemma \ref{9} that for all $j\in \{0,\ldots,k\}$, the series $\sum_m \hat{\phi}_m \mathrm{ad}_A^j U^m$ is norm convergent. The conclusion follows by applying Lemma \ref{H0-2}. \ep
\\

\noindent{\bf Remark:} Propositions \ref{4}, \ref{10} and Corollary \ref{8-2} allows us to derive an alternative proof of \cite{abcf} Lemma 5.2.

\noindent{\bf Remark:} It is also clear that for any bounded symmetric operator $B$ defined on ${\cal H}$, any $\alpha \in {\mathbb R}\setminus \{0\}$ and any $k\in {\mathbb N}$, $C^k(A)=C^k(\alpha A+B)$.

\subsection{The class ${\cal C}^{1,1}(A)$ and approximation properties}

If $B$ is a bounded linear operator on ${\cal H}$, which belongs to $C^1(A)$, we denote by ${\cal A}_0(B)$ the class of families $(B_{\epsilon})_{\epsilon \in (0,\epsilon_0]}$ (for some $\epsilon_0 >0$) such that:
\begin{itemize}
\item[(a)] There exists $c > 0$ such that $\|B_{\epsilon}-B\|\leq c \, \epsilon$,
\item[(b)] the function $\epsilon \to B_{\epsilon}$ belongs to $C^1((0,\epsilon_0], \mathcal{B}(\mathcal{H}))$,
\item[(c)] the map $\epsilon \mapsto \epsilon^{-1}\|\partial_{\epsilon}B_{\epsilon}\|$ belongs to $L^1([0,\epsilon_0])$.
\end{itemize}
The class ${\cal A}_1(A,B)$ will denote the subclass of families $(B_{\epsilon})_{\epsilon \in (0,\epsilon_0]}\in {\cal A}_0(B)$ such that:
\begin{itemize}
\item[(a)] for all $\epsilon \in (0,\epsilon_0]$, $B_{\epsilon} \in C^1(A)$,
\item[(b)] $\|[A,B_{\epsilon}]-[A,B]\|$ vanishes as $\epsilon$ tends to $0$,
\item[(c)] The map $\epsilon \to [A,B_{\epsilon}]$ belongs to $C^1((0,\epsilon_0], \mathcal{B}(\mathcal{H}))$,
\item[(d)] $\|\partial_{\epsilon} [A,B_{\epsilon}] \|$ belongs to $L^1([0,\epsilon_0])$,
\end{itemize}
The class ${\cal A}_2(A,B)$ is the subclass of families $(B_{\epsilon})_{\epsilon \in (0,\epsilon_0]}\in {\cal A}_1(A,B)$ such that:
\begin{itemize}
\item[(a)] for all $\epsilon \in (0,\epsilon_0]$, $B_{\epsilon} \in C^2(A)$
\item[(b)] $\|\mathrm{ad}_A^2 B_{\epsilon}\|$ belongs to $L^1([0,\epsilon_0])$.
\end{itemize}

This leads us to the following definition:
\begin{defin} Let $B$ is a bounded linear operator on ${\cal H}$, which belongs to $C^1(A)$. $B\in {\cal C}(A)$ if ${\cal A}_2(A,B)\neq \emptyset$.
\end{defin}

The following proposition is immediate:
\begin{prop}\label{alg} Let $B$ and $C$ two operators which belong to ${\cal C}(A)$. Then for any $\lambda\in {\mathbb C},$
\begin{itemize}
\item $BC \in {\cal C}(A)$
\item $B+\lambda C \in {\cal C}(A)$
\item $B^* \in {\cal C}(A)$
\end{itemize}
\end{prop}
\noindent{\bf Proof:} It is enough to see that if $B$ and $C$ are operators which belong to $C^1(A)$ and $(B_{\epsilon})_{\epsilon \in (0,\epsilon_0]} \in {\cal A}_k(A,B)$ and $(C_{\epsilon})_{\epsilon \in (0,\epsilon_1]}\in {\cal A}_k(A,C)$ ($k\in \{0,1,2\}$), then for any $\lambda\in {\mathbb C},$
\begin{itemize}
\item $(B_{\epsilon}C_{\epsilon})_{\epsilon \in (0,\min(\epsilon_0,\epsilon_1)]} \in {\cal A}_k(A,BC)$
\item $(B_{\epsilon}+\lambda C_{\epsilon})_{\epsilon \in (0,\min(\epsilon_0,\epsilon_1)]} \in {\cal A}_k(A,B+\lambda C)$
\item $(B_{\epsilon}^*)_{\epsilon \in (0,\epsilon_0]} \in {\cal A}_k(A,B^*)$.
\end{itemize}
\ep

The central point of this discussion lies in the following result:
\begin{prop}\label{georgescu} ${\cal C}^{1,1}(A)= {\cal C}(A)$.
\end{prop}
\noindent{\bf Proof:} If $B$ is any bounded symmetric operator, $B\in {\cal C}^{1,1}(A)$ if and ony if $B\in {\cal C}(A)$. This follows from \cite{abmg} paragraph 7.3 and Lemma 7.3.6. Since ${\cal C}^{1,1}(A)$ and ${\cal C}(A)$ are stable under $*$, the result follows by considering the real and imaginary parts and using the vector space structure of these classes. \ep

However, for a given unitary operator, we would like to express this approximation property in terms of unitary operators as it appears in the proof of Theorem \ref{thm0}. Let us recall first the following result:
\begin{prop}\label{hscor} Let $B$ a bounded self-adjoint operator defined on ${\cal H}$ and $\Phi \in C_0^{\infty}({\mathbb R})$. If $B\in {\cal C}^{1,1}(A)$ then $\Phi(B)\in {\cal C}^{1,1}(A)$.
\end{prop}
See \cite{abmg} Theorem 6.2.5, Corollary 6.2.6 or \cite{gj} for a proof.

Now, we show that if a unitary operator $U$ belongs ${\cal C}^{1,1}(A)$, then the approximating family can be chosen to be unitary.
\begin{prop}\label{u} Let $U$ a unitary operator which belongs to $C^1(A)$ and $k\in \{1,2\}$. If $(B_{\epsilon})_{\epsilon \in (0,\epsilon_0]} \in {\cal A}_0(U)$ (resp. ${\cal A}_k(A,U)$), then there exists a family of unitary operators $(U_{\epsilon})_{\epsilon \in (0,\epsilon_1]} \in {\cal A}_0(U)$ (resp. ${\cal A}_k(A,U)$) for some $\epsilon_1 \in (0,\epsilon_0]$.
\end{prop}
\noindent{\bf Proof:} Let us consider the (unique) polar decomposition of $B_{\epsilon}$:
\begin{equation*}
B_{\epsilon}= |B_{\epsilon}| U_{\epsilon} \qquad \mathrm{where}\quad |B_{\epsilon}| = \sqrt{B_{\epsilon}^*B_{\epsilon}}
\end{equation*}
and $(U_{\epsilon})$ is a family of partial isometry \cite{Kato1}. Since $U$ is unitary and
\begin{equation*}
\|B_{\epsilon} -U\| \leq C \epsilon \enspace ,
\end{equation*}
there exists $\epsilon_1 \in (0,\epsilon_0]$ such that for all $\epsilon \in (0,\epsilon_1]$, $B_{\epsilon}$ is also invertible, which implies that $(U_{\epsilon})_{\epsilon \in (0,\epsilon_1]}$ is a family of unitary operators that can be defined (in the norm sense) by:
\begin{eqnarray*}
U_{\epsilon} &=& |B_{\epsilon}|^{-1} B_{\epsilon}\\
|B_{\epsilon}|^{-1} &=& \sqrt{B_{\epsilon}^* B_{\epsilon}}^{-1}=\Phi(B_{\epsilon}^* B_{\epsilon})
\end{eqnarray*}
for a suitable $\Phi \in C^{\infty}_0({\mathbb R})$. By Proposition \ref{alg}, the family $(B_{\epsilon}^* B_{\epsilon})$ also belongs to ${\cal A}_0(U)$ (resp. ${\cal A}_k(A,U)$). Due to Lemma \ref{hscor}, this also the case for the family $(\Phi(B_{\epsilon}^* B_{\epsilon}))$. Applying again Proposition \ref{alg}, we have that $(U_{\epsilon})_{\epsilon \in (0,\epsilon_1]}$ also belongs to ${\cal A}_0(U)$ (resp. ${\cal A}_k(A,U)$). \ep

Rewriting Propositions \ref{georgescu} and \ref{u}, we obtain that:
\begin{cor}\label{u2} Let $U$ be a unitary operator on ${\cal H}$. Then, $U$ belongs to ${\cal C}^{1,1}(A)={\cal C}(A)$ if and only if there exists a family of unitary operators $(U_{\epsilon})_{\epsilon \in (0,\epsilon_0]}$ (for some $\epsilon_0 >0$) such that:
\begin{itemize}
\item[(a)] for all $\epsilon \in (0,\epsilon_0]$, $U_{\epsilon}\in C^2(A)$,
\item[(b)] there exists $c > 0$ such that $\|U_{\epsilon}-U\|\leq c\epsilon$,
\item[(c)]  $\|[A,B_{\epsilon}]-[A,B]\|$ vanishes as $\epsilon$ tends to $0$,
\item[(d)] the map $\epsilon \to U_{\epsilon}$ is $C^1$ on $(0,\epsilon_0]$ with respect to the norm topology on $\mathcal{B}(\mathcal{H})$,
\item[(e)] the map $\epsilon \to [A,U_{\epsilon}]$ is $C^1$ on $(0,\epsilon_0]$ with respect to the norm topology on $\mathcal{B}(\mathcal{H})$,
\item[(f)] the map $\epsilon \mapsto \epsilon^{-1}\|\partial_{\epsilon}U_{\epsilon}\|+\| \partial_{\epsilon}[A,U_{\epsilon}]\|+\|\mathrm{ad}_A^2 B_{\epsilon}\|$ belongs to $L^1([0,\epsilon_0])$.
\end{itemize}
\end{cor}

\subsection{On the proof of Lemma \ref{lemEnss3}}

We have chosen to prove this lemma by using the correspondence described in Lemma \ref{georgescu} and Corollary \ref{u2}. The context of this discussion is described in Section 3 to which we refer. We recall that for $V(\cdot) \in C^k({\mathbb R})$, the multiplication operator by $\partial_x^k V(\cdot)$ is denoted by $\partial_x^k V$ or equivalently $\mathrm{ad}_p^k V$.

If $V\in {\cal C}^{1,1}(p)\subset C^1(p)$, we know by \cite{abmg} Lemma 7.3.6, there exists a family of real-valued functions $(V_{\epsilon}(\cdot))_{\epsilon \in (0,1]}$ such that the associated family of multiplication operators on $L^2({\mathbb R})$, $(V_{\epsilon}(\cdot))_{\epsilon \in (0,1]}$, belongs to ${\cal A}_2(p,V)$. Namely,
\begin{equation*}
V_{\epsilon}(x) = \int_{-\infty}^{\infty} V(x-\epsilon \tau) e^{-\frac{\tau^2}{4}}\, \frac{d\tau}{\sqrt{4\pi}} \enspace .
\end{equation*}
It is natural to define the unitary propagator associated to $V_{\epsilon}$ by $U_{\epsilon}(t,s) = U_0(t,s) \Omega_{\epsilon}(t,s)$ for all $(s,t)\in {\mathbb R}^2$, where $\Omega_{\epsilon} (t,s)$ is defined in the strong sense by:
\begin{equation*}
\Omega_{\epsilon}(t,s)-I = -i \int_s^t U_0(\tau,s) V_{\epsilon} U_0^*(\tau,s) \Omega_{\epsilon}(\tau,s)\, d\tau \enspace.
\end{equation*}
In addition, according to Lemma \ref{lemEnss1}, since $V(\cdot)$ and $V_{\epsilon}(\cdot)$ belong to $C^1(p)$ we can define the bounded commutators:
\begin{eqnarray*}
B_{\epsilon}(T) &=& U_{\epsilon}^*(T)pU_{\epsilon}(T)-p\\
B(T) &=& U^*(T)pU(T)-p \enspace,
\end{eqnarray*}
with the conventions: $U(T):=U(T,0)$ and $U_{\epsilon}(T):=U_{\epsilon}(T,0)$.

As a consequence we obtain that:
\begin{lem}\label{approx1} Assume there exists a family of real-valued function $(V_{\epsilon}(\cdot))_{\epsilon \in (0,1]}$ such that $(V_{\epsilon})_{\epsilon \in (0,1]} \in {\cal A}_0(V)$. Then, for each $(s,t)\in {\mathbb R}^2$, $(U_{\epsilon}(t,s))_{\epsilon \in (0,1]} \in {\cal A}_0(U(t,s))$.
\end{lem}
\noindent{\bf Proof:} Fix $(s,t)\in {\mathbb R}^2$. By Duhamel formula, we have that 
\begin{equation}\label{duhamel}
U_{\epsilon}(t,s) - U(t,s) = i\int_s^t U_{\epsilon}(t,\tau)(V_{\epsilon} - V) U(\tau,s) \,d\tau \enspace.
\end{equation}
The hypotheses and the unitarity of the propagators $(U_{\epsilon}(t,s))$ and $(U(t,s))$ imply that ${\cal A}_0(U(t,s))$ Statement (a) holds for $U_{\epsilon}(t,s)$. Duhamel formula also rewrites:
\begin{equation*}
U_{\epsilon_1}(t,s) - U_{\epsilon_2}(t,s) = i\int_s^t U_{\epsilon_1}(t,\tau)(V_{\epsilon_1} - V_{\epsilon_2}) U_{\epsilon_2}(\tau,s) \,d\tau \enspace,
\end{equation*}
for any $(\epsilon_1,\epsilon_2)\in (0,1]^2$, which once combined with the hypotheses shows that the map $\epsilon \mapsto U_{\epsilon}(t,s)$ is norm-continuous. On the same basis, this map is differentiable and for any $\epsilon \in (0,1)$ and
\begin{equation*}
\partial_{\epsilon} U_{\epsilon}(t,s) = i \int_s^t U_{\epsilon}(t,\tau) \left( \partial_{\varepsilon} V_{\epsilon}\right) U_{\epsilon}(\tau,s) \,d\tau \enspace .
\end{equation*}
${\cal A}_0(U(t,s))$ Statement (b) follows as above. Combining the last identity with the hypotheses also entails that:
\begin{equation*}
\int_0^1 \epsilon^{-1} \| \partial_{\epsilon} U_{\epsilon}(t,s) \|  \,d\epsilon \leq |t-s| \int_0^1 \epsilon^{-1} \| \partial_{\epsilon}V_{\epsilon} \|  \,d\epsilon \enspace ,
\end{equation*}
which proves ${\cal A}_0(U(t,s))$ Statement (c). \ep

\begin{lem}\label{approx2} Assume there exists a family of real-valued function $(V_{\epsilon}(\cdot))_{\epsilon \in (0,1]}$ such that $(V_{\epsilon})_{\epsilon \in (0,1]}$ belongs to ${\cal A}_1(p,V)$ (resp. ${\cal A}_2(p,V)$). Then, $(U_{\epsilon}(T))_{\epsilon \in (0,1]}$ belongs to ${\cal A}_1(p,U(T))$ (resp. ${\cal A}_2(p,U(T))$).
\end{lem}
\noindent{\bf Proof:} We deduce from Lemma \ref{approx1} that $(U_{\epsilon}(T))_{\epsilon \in (0,1]}$ belongs to ${\cal A}_0(U(T))$. Note also that Statements ${\cal A}_1(p,U(T))$ (a) and ${\cal A}_2(p,U(T))$ (a) follow from Lemma \ref{lemEnss1}. In particular, we have that:
\begin{equation*}
B_{\epsilon}(T) -B(T) = \int_0^T \cos(\omega\tau) \left(U^*(\tau) \left(\partial_x V\right) U(\tau) - U_{\epsilon}^*(\tau) \left(\partial_x V_{\epsilon}\right) U_{\epsilon}(\tau)\right) \, d\tau \enspace ,
\end{equation*}
for all $\epsilon \in (0,1]$. Statement ${\cal A}_1(p,U(T))$ (b) is immediate. We also have that:
\begin{equation*}
B_{\epsilon_1}(T) -B_{\epsilon_2}(T) = \int_0^T \cos(\omega\tau) \left(U_{\epsilon_2}^*(\tau) \partial_x V_{\epsilon_2} U_{\epsilon_2}(\tau) - U_{\epsilon_1}^*(\tau) \partial_x V_{\epsilon_1} U_{\epsilon_1}(\tau)\right) \, d\tau \enspace ,
\end{equation*}
for any $(\epsilon_1,\epsilon_2)\in (0,1]^2$. A similar argument to the proof of Lemma \ref{approx1} allows to prove Statement ${\cal A}_1(p,U(T))$ (b) together with:
\begin{equation*}
\partial_{\varepsilon} B_{\epsilon}(T) = -\int_0^T \cos(\omega\tau) U_{\epsilon}^*(\tau) \left(\partial_{\varepsilon}(\partial_x V_{\epsilon})\right) U_{\epsilon}(\tau)\, d\tau \enspace,
\end{equation*}
for all $\epsilon \in (0,1]$. Then, Statement ${\cal A}_1(p,U(T))$ (d) follows directly. We end up with the lemma rewriting $B_{\epsilon}(T)$ as:
\begin{equation*}
B_{\epsilon}(T) = -\int_0^T \cos(\omega\tau) \Omega_{\epsilon}^*(\tau) U_0^*(\tau) V_{\epsilon}^{(1)}(\tau)U_0(\tau)\Omega_{\epsilon}(\tau) \, d\tau \enspace .
\end{equation*}
It follows from the hypotheses and Lemma \ref{lemEnss0-1} that $\|\mathrm{ad}_p \Omega_{\epsilon}^*(\tau) \|$ , $\| \mathrm{ad}_p \Omega_{\epsilon}(\tau) \|$
are bounded by $T \|\partial_x V_{\epsilon}\|$. it also follows from Lemma \ref{lemEnss0} that:
$\|\mathrm{ad}_p (U_0^*(\tau) \partial_x V_{\epsilon}(\tau) U_0(\tau)) \| \leq  \|\partial_x^2 V_{\epsilon}\|$ on the interval $0 \leq \tau \leq T$. This allows to prove Statement ${\cal A}_2(p,U(T))$ (b): by hypotheses one has that $\|\partial_x V_{\epsilon}\|$ and $\|\partial_x^2 V_{\epsilon}\|$ are integrable with respect to $\epsilon$, $\epsilon \in (0,1]$ , thus $\|\mathrm{ad}_p (B_{\epsilon}(T)) \|$ is integrable in [0,1]. \ep
\\

\noindent{\bf Proof of Lemma \ref{lemEnss3}:} This is a combination of Proposition \ref{georgescu} and Lemma \ref{approx2}. \ep

\section{On Weyl Commutation Relations}

Let ${\cal H}$ be a Hilbert space. Let $T$ be a unitary operator defined on ${\cal H}$ and $A$ a self-adjoint operator with domain ${\cal D}(A)$ defined on ${\cal H}$. We denote by ${\cal G}_T$ (resp. ${\cal G}_A$) the discrete group generated by $T$ (resp. strongly continuous group generated by $A$), i.e.:
\begin{eqnarray*}
{\cal G}_T &=& \bra T^n; n\in {\mathbb Z} \ket\\
{\cal G}_A &=& \bra e^{itA}; t\in {\mathbb R} \ket
\end{eqnarray*}

In this context, we adopt the following definition:
\begin{defin} The groups ${\cal G}_T$ and ${\cal G}_A$ constitutes a Heisenberg couple if they satisfy the following Heisenberg commutation relation: for all $(t,k) \in {\mathbb R}\times {\mathbb Z}$, $e^{itA} T^k = e^{itk} T^k e^{itA}$.
\end{defin}

The following equivalence is straightforward:
\begin{prop}\label{Heisenberg} Let ${\cal H}$ be a Hilbert space. Let $T$ be a unitary operator defined on ${\cal H}$ and $A$ a self-adjoint operator defined on ${\cal H}$ with domain ${\cal D}(A)$. Then, the following assertions are equivalent:
\begin{itemize}
\item The sesquilinear form defined on ${\cal D}(A)\times {\cal D}(A)$ by $G(\varphi,\phi):=\langle T\varphi,AT\phi\rangle\ - \langle\varphi,A\phi\rangle$ extends continuously to a bounded form on ${\cal H}\times {\cal H}$. The associated bounded operator is $T^{*}AT-A=I$.
\item $T\in C^1(A)$ and $\mathrm{ad}_A T=[A,T]=T$.
\item For all $t\in {\mathbb R}$, $F(t) \equiv e^{itA} T e^{-itA} = e^{it}T$.
\item The groups ${\cal G}_T$ and ${\cal G}_A$ constitute an Heisenberg couple on ${\cal H}$.
\end{itemize}
\end{prop}
The equivalence between the two first statements has been established previously. Actually, the first statement is also equivalent to: 
\begin{equation*}
T^{-n} AT^n -A = n I \enspace,
\end{equation*}
for all $n\in {\mathbb Z}$. In other words $A$ is a time-operator for $T$, or more exactly for ${\cal G}_T$. This concept of time-operator is not limited to the Hilbert space setting and can be defined for discrete and strongly time-continuous semi-groups on Banach spaces \cite{such1}. \\

Actually, the function $F$ defined in Proposition \ref{Heisenberg} is norm analytic. This implies in particular that: $T\in C^{\infty}(A)$. It also follows from the same identity that for all $t\in {\mathbb R}$: $\sigma (T) = \sigma (e^{itA} T e^{-itA}) = e^{it} \sigma(T)$. So, $\sigma (T)=\partial {\mathbb D}={\mathbb S}^1$.

\section{$U$-smoothness and beyond}

We start by recalling basics facts about the concept of $U$-smoothness. Let $U$ be a unitary operator acting a fixed Hilbert space ${\cal H}$ and ${\cal F}$ an auxiliary Hilbert space. Then,
\begin{defin}\label{def-1} Let $B \in {\cal B}({\cal H},{\cal F})$ be a bounded operator. The operator $B$ is $U$-smooth if:
\begin{equation*}
sup_{\| \psi\|=1} \sum_{n\in {\mathbb Z}} \|BU^n\psi\|^2 < \infty\enspace.
\end{equation*}
\end{defin}

The notion of $U$-smoothness may also be defined locally. If $E_U(\cdot)$ denotes the spectral family associated to the unitary operator $U$, then given a Borel subset $\Theta$ of the one-dimensional torus ${\mathbb T}$, a bounded operator $B\in {\cal B}({\cal H},{\cal F})$, is locally $U$-smooth on $\Theta$ if the operator $BE_U(\Theta)$ is $U$-smooth. In particular, we say that $B$ is $U$-smooth at a point $\theta$ in ${\mathbb T}$, if there exists an open neigbourhood $\Theta_{\theta}$ of $\theta$ such that $BE_U(\Theta_{\theta})$ is $U$-smooth. For a proof of the following equivalence, see \cite{abcf}:
\begin{thm}\label{characterization} Let $U$ be a unitary operator and $B$ be a bounded operator defined on the Hilbert space ${\cal H}$. If the operator $B$ is $U$-smooth, then the statements
\begin{enumerate}
\item \quad $$C_1 \equiv \frac{1}{2\pi}\sup_{\|\psi\|=1} \sum_{n\in {\mathbb Z}} \| BU^n\psi\|^2 < \infty$$
\item \quad $$C_2 \equiv \frac{1}{2\pi}\sup_{\|\psi\|=1, z\in {\mathbb D}} |\langle B^{\ast}\psi ,\Re (\frac{1+zU^*}{1-zU^*}) B^{\ast}\psi \rangle | < \infty$$
\item \quad $$C_3 \equiv \sup_{\substack{\|\psi\|=1\\(a,b)\in {\mathbb T}^2,a<b}}\dfrac{\| BE_U((a,b))\psi\|^2}{|b-a|}<\infty$$
\item \quad $$C_4 \equiv \sup_{\substack{\|\psi\|=1\\(a,b)\in {\mathbb T}^2,a<b}}\dfrac{\| E_U((a,b))B^{\ast}\psi\|^2}{|b-a|}<\infty$$
\item \quad $$C_5 \equiv \frac{1}{2\pi} \sup_{\|\psi\|=1,z\in {\mathbb D}}(1-|z|^2) \|(1-zU^*)^{-1}B^{\ast}\psi\|^2 <\infty$$
\end{enumerate}
are equivalent. Moreover: $C_1=C_2=C_3=C_4=C_5$.
\end{thm}

We proceed now to some extensions of this concept. The following discussion is adapted from \cite{abmg} Paragraph 7.1.

Let ${\cal K}$ be a Banach space, continuously and densely embedded in ${\cal H}$. Identifying the Hilbert space ${\cal H}$ with its dual, we have that: ${\cal H}\subset {\cal K}^*$. Denote by ${\cal K}'$ the closure of ${\cal H}$ in ${\cal K}^*$ equipped with its Banach space topology. So, ${\cal H}$ is continuously and densely embedded in ${\cal K}'$ and ${\cal K}'$ is a closed subspace of ${\cal K}^*$. It follows that ${\cal B}({\cal H})$ is embedded in ${\cal B}({\cal K},{\cal K}')$, which is itself isometrically embedded in ${\cal B}({\cal K},{\cal K}^*)$. ${\cal B}({\cal K},{\cal K}')$ is a norm-closed, weak$^*$ dense subspace of ${\cal B}({\cal K},{\cal K}^*)$. This means in particular  that the operators defined for $|z|<1$ by $(1-zU^*)^{-1}$, $(1-\bar{z}^{-1}U^*)^{-1}$ and
\begin{equation*}
\Re (\frac{1+zU^*}{1-zU^*})= (1-zU^*)^{-1} - (1-\bar{z}^{-1}U^*)^{-1}
\end{equation*}
belong to ${\cal B}({\cal H})$ but also to ${\cal B}({\cal K},{\cal K}^*)$.

Let $\Theta$ be an open subset of ${\mathbb T}$ and define the following properties:
\begin{itemize}
\item[(U1)] For any compact subset $K\subset \Theta$, there exists $C_K>0$, such that for any $\psi \in {\cal K}$ and any $|z|<1$ with $\arg z \in K,$
\begin{equation*}
|\bra \psi, \Re (\frac{1+zU^*}{1-zU^*}) \psi \ket | \leq C_K \|\psi\|^2_{\cal K} \enspace .
\end{equation*}
\item[(U2)] For all $\theta \in \Theta$,
\begin{equation*}
\lim_{z\rightarrow e^{i\theta}, |z|<1} \Re (\frac{1+zU^*}{1-zU^*})
\end{equation*}
exists (uniformly in $\theta$ on each compact subset $K\subset \Theta$) in the weak$^*$ topology of ${\cal B}({\cal K},{\cal K}^*)$.
\end{itemize}
\noindent{\bf Remark:} Property (U1) is equivalent to: for all compact subset $K\in \Theta,$
\begin{equation*}
\sup_{\arg z\in K, |z|<1}\|\Re (\frac{1+zU^*}{1-zU^*})\|_{{\cal B}({\cal K},{\cal K}^*)} < \infty \enspace .
\end{equation*}
The Banach-Steinhaus Theorem shows that property (U2) implies property (U1). The next proposition is a straightforward adaptation of \cite{abmg} Lemma 7.1.3 (see also the proof of Theorem \ref{characterization} in \cite{abcf}).
\begin{prop}\label{usmoothext} If ${\cal K}$ is continuously and densely embedded in ${\cal H}$, we have that:
\begin{itemize}
\item If property (U1) holds, then $U$ is purely absolutely continuous on $\Theta$.
\item If property (U2) holds, then for all $\theta_0 \in {\mathbb T}$ fixed, the map $\theta \mapsto E_{\theta} - E_{\theta_0}$ is differentiable with respect to the weak $^*$ topology on ${\cal B}({\cal K},{\cal K}^*)$ and
\begin{equation*}
\frac{dE_{\theta}}{d\theta} = \frac{1}{2\pi}\lim_{r\mapsto 1^-} \Re (\frac{1+re^{i\theta}U^*}{1-re^{i\theta}U^*}) = \frac{1}{2\pi}\lim_{r\mapsto 1^-} (1-re^{i\theta}U^*)^{-1} - (1-r^{-1}e^{i\theta}U^*)^{-1}\enspace .
\end{equation*}
\item Assume property (U1) holds and $({\cal K}')^*= {\cal K}$. If $T\in {\cal B}({\cal K}',{\cal F})$ then $T$ is locally $U$-smooth on $\Theta$.
\end{itemize}
\end{prop}

\noindent{\bf Acknowledgments:} The authors thank Joachim Asch for informative discussions.


\end{document}